\documentclass[12pt]{iopart}
\expandafter\let\csname equation*\endcsname\undefined
\expandafter\let\csname endequation*\endcsname\undefined
\usepackage{geometry}
\usepackage{amsthm}
\usepackage{dsfont}
\usepackage{bm, upgreek}

\usepackage{placeins}

\usepackage{amsmath}
\usepackage[title]{appendix}
\usepackage{graphicx}
\usepackage{subcaption}
\usepackage{iopams}
\usepackage{esint}

\usepackage{siunitx}
\usepackage{algorithm}
\usepackage{algpseudocode}
\usepackage{booktabs}

\usepackage[outdir=./Images/]{epstopdf}
\usepackage{wrapfig}

\usepackage[backend=biber,sortcites]{biblatex}
\bibliography{biblio}

\newtheorem{theorem}{Theorem}
\newtheorem{lemma}{Lemma}
\newtheorem{cor}{Corollary}

\theoremstyle{remark}
\newtheorem*{remark}{Remark}

\DeclareMathOperator*{\argmin}{argmin}
\DeclareMathOperator*{\diag}{diag}

\DeclareMathOperator*{\minimise}{minimise}

\DeclareMathOperator{\prox}{prox}
\DeclareMathOperator{\dom}{dom}
\DeclareMathOperator{\proj}{proj}
\DeclareMathOperator{\Span}{span}

\newcommand{\TGV}{\mathrm{TGV}}
\newcommand{\TV}{\mathrm{TV}}
\newcommand{\VTV}{\mathrm{VTV}}
\newcommand{\imi}{\mathrm{i}}

\newcommand*{\lvec}[1]{\boldsymbol{\mathbf{#1}}}
\newcommand*{\lmat}[1]{\mathrm{#1}}
\newcommand*{\lele}[1]{{#1}}

\newcommand*{\svec}[1]{\boldsymbol{#1}}

\usepackage{xargs}
\usepackage[pdftex,dvipsnames]{xcolor}
\usepackage[colorinlistoftodos,prependcaption,textsize=tiny]{todonotes}
\newcommandx{\unsure}[2][1=]{\todo[linecolor=Plum,backgroundcolor=Plum!25,bordercolor=Plum,#1]{#2}}
\newcommandx{\change}[2][1=]{\todo[linecolor=blue,backgroundcolor=blue!25,bordercolor=blue,#1]{#2}}
\newcommandx{\info}[2][1=]{\todo[linecolor=OliveGreen,backgroundcolor=OliveGreen!25,bordercolor=OliveGreen,#1]{#2}}

\begin{document}
	\title{On Regularisation of Coherent Imagery with Proximal Methods}
\author{FM Watson, J Hellier and WRB Lionheart}

\address{Department of Mathematics, University of Mancheter, Oxford Road, Manchester, M13 9PL, UK}
\ead{francis.watson@manchester.ac.uk}
\vspace{10pt}
\begin{indented}
	\item[]
\end{indented}

\begin{abstract}
	In complex-valued coherent inverse problems such as synthetic aperture radar (SAR), one may often have prior information only on the magnitude image which shows the features of interest such as strength of reflectivity.  In contrast, there may be no more prior knowledge of the phase beyond it being a uniform random variable.  However, separately regularising the magnitude, via some function \(G:=H(|\cdot|)\), would appear to lead to a potentially challenging non-linear phase fitting problem in each iteration of even a linear least-squares reconstruction problem.  We show that under certain sufficient conditions the proximal map of such a function \(G\) may be calculated as a simple phase correction to that of \(H\). Further, we provide proximal map of (almost) arbitrary \(G:=H(|\cdot|)\) which does not meet these sufficient conditions.  This may be calculated through a simple numerical scheme making use of the proximal map of \(H\) itself, and thus we provide a means to apply practically arbitrary regularisation functions to the magnitude when solving coherent reconstruction problems via proximal optimisation algorithms.  This is demonstrated using publicly available real SAR data for generalised Tikhonov regularisation applied to multi-channel SAR, and both a simple level set formulation and total generalised variation applied to the standard single-channel case.
\end{abstract}


\section{Introduction}
We can think of coherent image formation and inverse problems as complex-valued ones in which phase is accurately preserved. In such inverse and imaging problems it may often be desirable to apply regularisation to the magnitude of the image only; prior knowledge may indeed only be available for the structure of the magnitude, possibly with nothing more to say about the phase than it is a uniform random variable.  

In the case of synthetic aperture radar (SAR), the phase of the image is highly sensitive to small positional errors in a way in which the magnitude is not. As a result, this phase of a single SAR image is generally effectively meaningless when viewed in isolation, appearing largely random. Such positional errors include small offsets of a potentially gridded pixel location from the true scattering point, if indeed the interaction is well represented by being from a single point scatterer.  In general one will most often view the magnitude image alone (referred to as the ``detected'' image), which describes the reflectivity strength of the observed targets. 

Yet, phase differences between SAR image pairs provide very valuable information through interferometry, such as terrain height, the detection of very small changes between collections, or resolving target motion\cite{gens1996review, rocca2000sar, jakowatz2012spotlight, andre2024moving}.  It may therefore not only be the case that prior information is only available about the reflectivity (amplitude), but that we want to ensure regularisation methods do not affect the phase to preserve these phase relationships, or even that a (poor) prior on phase may prevent a good reconstruction at all. 

\begin{figure}[h!]
	\centering
	\begin{subfigure}[t]{0.32\textwidth}
		\centering
		\includegraphics[width=0.99\textwidth, trim=2cm 3cm 3cm 3cm, clip]{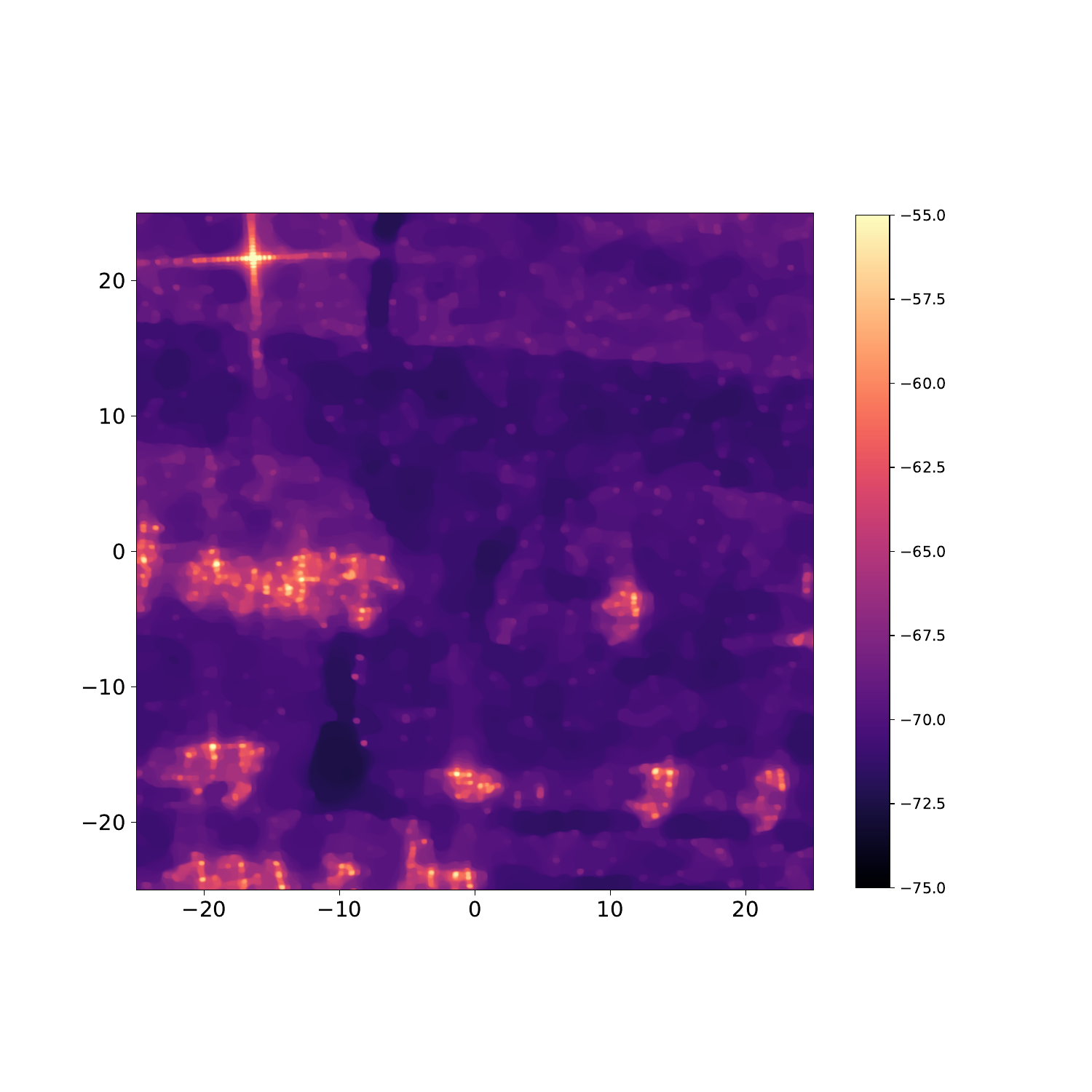}
		\caption{\(\TV(|\cdot|)\)-regularised}
		\label{subfig: TVabs reg gotcha}
	\end{subfigure}
	\begin{subfigure}[t]{0.32\textwidth}
		\centering
		\includegraphics[width=0.99\textwidth, trim=2cm 3cm 3cm 3cm, clip]{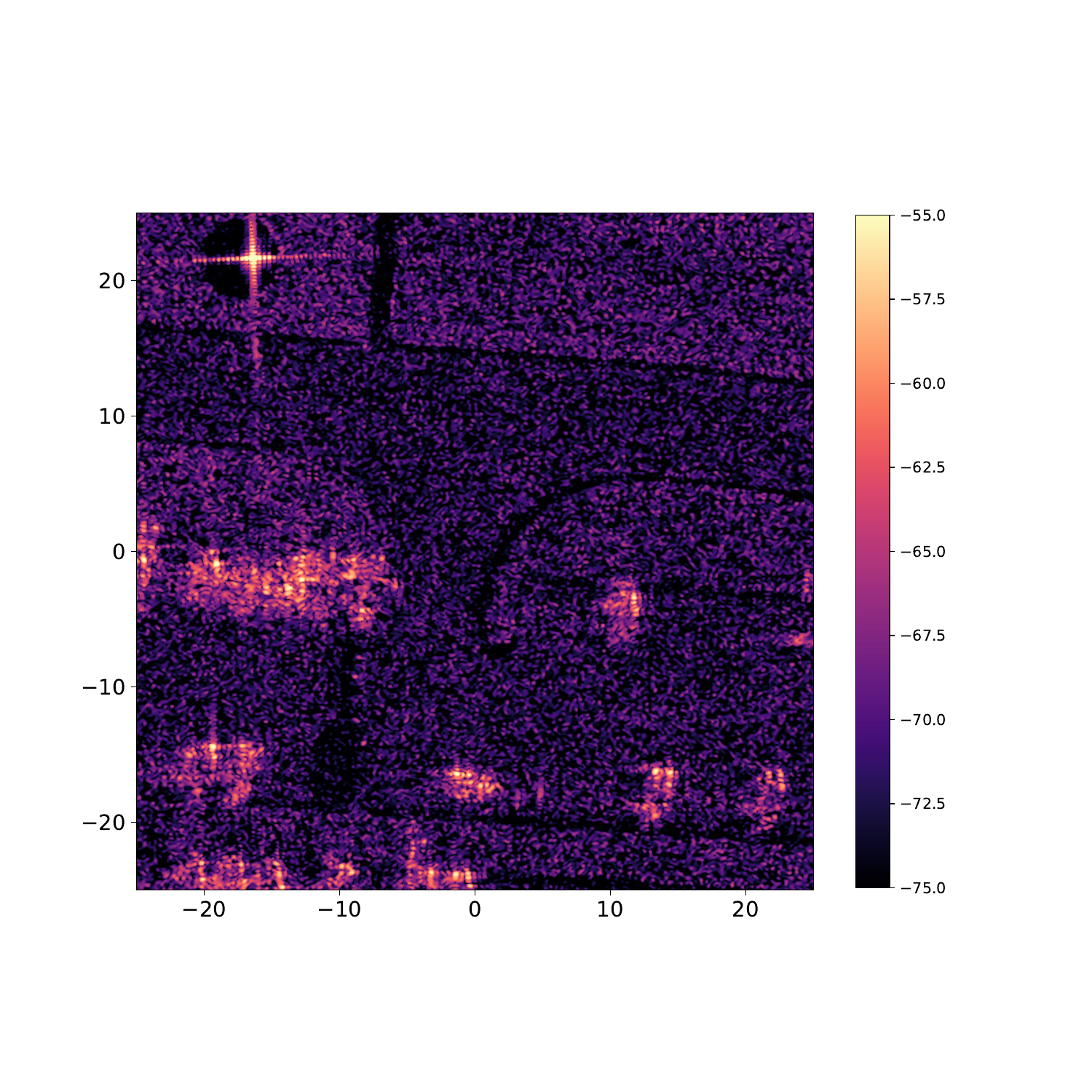}
		\caption{\(\TV\)-regularised}
		\label{subfig: TV reg gotcha}
	\end{subfigure}
	\begin{subfigure}[t]{0.32\textwidth}
		\centering
		\includegraphics[width=0.99\textwidth, trim=2cm 3cm 3cm 3cm, clip]{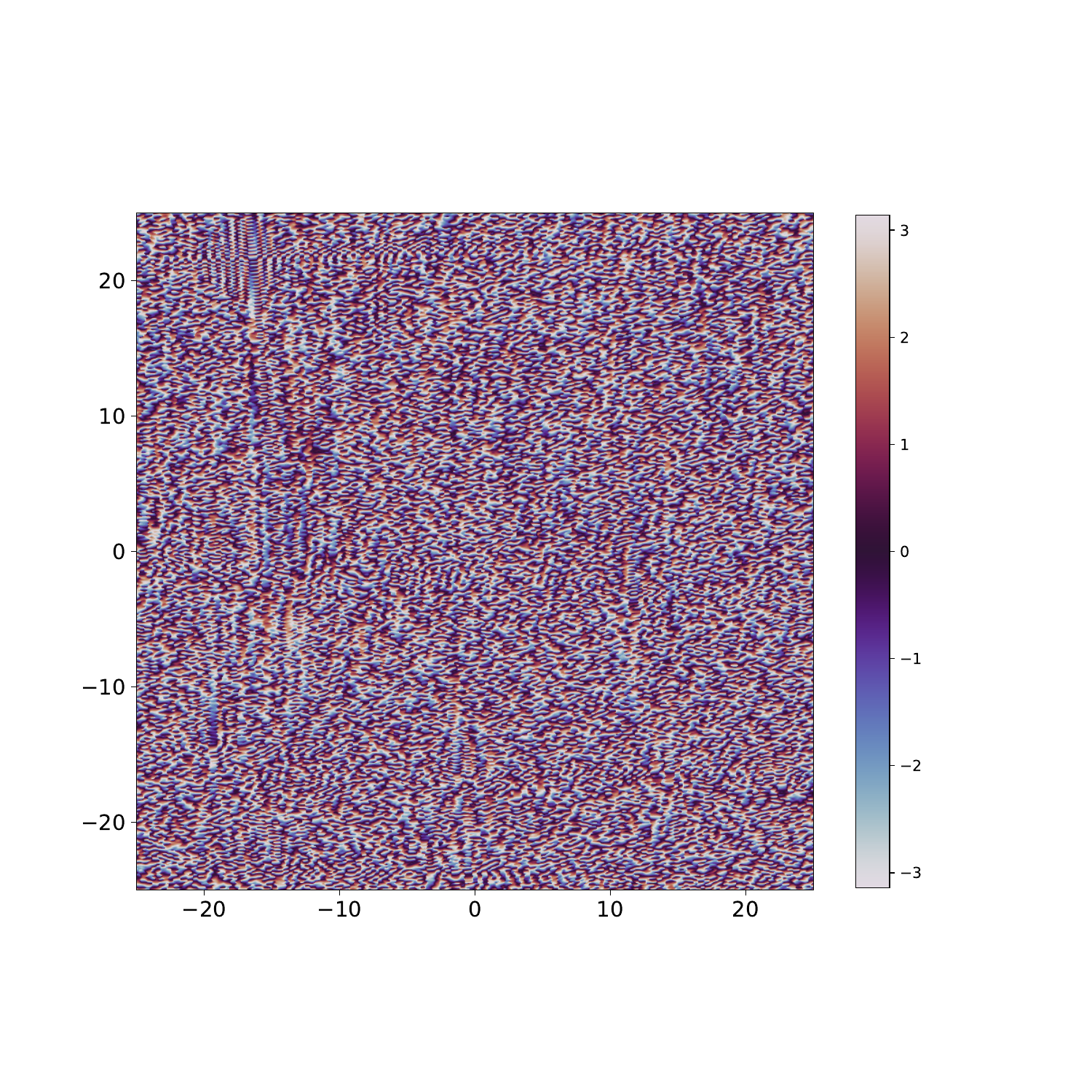}
		\caption{Phase of backprojection}
		\label{subfig: phase gotcha}
	\end{subfigure}
\caption{Example magnitude image chips of reconstructions of the Gotcha carpark data, showing (a) the \(\TV(|\cdot|)\) regularised result and (b) the result of applying \(\TV\) directly to the complex-valued image. Also shown in (c) is the phase of the backprojection image.}
\label{fig: tv into}
\end{figure}

As a motivating example, consider total variation (TV) regularised least-squares reconstruction, which we discuss in more detail in section~\ref{sec: common regs}. \figurename~\ref{fig: tv into} shows TV-regularised reconstruction image chips of the Gotcha carpark data\cite{casteel2007challenge}.  Applying TV to the magnitude directly, following the method later derived in this work as well as has been previously derived specifically for TV\cite{guven2016augmented}, we see from \figurename~\ref{subfig: TVabs reg gotcha} results in the expected and classic piecewise-constant appearance of a TV regularised reconstruction. This includes a reduction of speckle clutter. However, when applying TV to the complex values directly, for example as has been considered by Aghamiry \textit{et al} in the case of seismic imaging\cite{aghamiry2021complex}\footnote{They refer to this case as ``Treating \(x\) as a Real-Valued Variable''}, the large pixel-to-pixel phase changes required to fit the data means both a small data misfit and total variation together are not possible. There is seemingly no speckle reduction, nor a piecewise-constant appearance as the use of TV should promote. In fact, we found that both the data misfit and regularisation term were an order of magnitude greater at the minimum than in the case of applying TV to the magnitude: the two functions simply fight one another.

Thus, the motivation of this work is to provide a means to easily apply regularisation functions of the form \(G(\lvec{z}):=H(|\lvec{z}|)\), where \(|\cdot|\) is understood to mean element-wise absolute value, in complex-valued reconstruction problems. However, in parameterising the complex-valued model parameters \(\lvec{z}\) in terms of magnitude and phase, \(\lvec{z}=\lvec{r}\rme^{\imi \lvec{\phi}}\), this would seemingly lead to a nonlinear problem in recovering the phase \(\lvec{\phi}\) due to its appearance in an exponential.  This has previously been addressed, for example, by carrying out an alternating optimisation in \(\lvec{r}\) and \(\lvec{\phi}\) for applications in MRI\cite{zhao2012separate}, seismic imaging\cite{aghamiry2021complex}, and SAR\cite{rambour2019introducing}, in order to separate a simpler linear reconstruction problem in \(\lvec{r}\).  The application of TV to the magnitude SAR image formation was also originally carried out via a modified non-linear quasi-Newton method\cite{cetin2001feature}. More recently however, G{\"u}ven \textit{et al}\cite{guven2016augmented} have presented a method based on the efficient calculation of the proximal map of \(\TV(|\cdot|)\), \(\prox_{\TV(|\cdot|)}(\lvec{z})\), which avoids the nonlinear phase fitting problem.  This is given by a simple phase change to the map \(\prox_{\TV}(\lvec{r})\circ\rme^{\imi\lvec{\phi}}\).

Our main result is, firstly, to show that this same simple proximal map of \(G\) exists for many such functions \(G(\lvec{z}):= H(|\lvec{z}|)\) under certain sufficient conditions. This result, set out in section~\ref{sec: main result}, allows efficient primal-dual proximal optimisation algorithms to be applied, which solve problems of the form
\begin{equation}
	x_{\mathrm{sol}} = \argmin_x F(x) + G(x).
	\label{eq: optimisation general form}
\end{equation} 
Proximal optimisation algorithms such as FISTA\cite{beck2009siam, beck2009ieee}, PDHG\cite{chambolle2011first} and ADMM\cite{parikh2014proximal} can each be used to solve efficiently solve various problems with a convex function \(F\) and possibly non-smooth \(G\), provided the proximal map of \(G\) is itself readily available in closed form (i.e. it is \textit{proximable}) or it is otherwise easy to compute numerically.  
By providing a simple means to calculate the proximal of \(G(\lvec{z}):= H(|\lvec{z}|)\), we can solve the problem of form (\ref{eq: optimisation general form}) via proximal algorithms, we forgo the need to solve a non-linear phase fitting problem altogether.

With this first result, in section~\ref{sec: common regs} we demonstrate how some commonly used regularisation functions meet these sufficient conditions, and thus may easily be applied directly to the magnitude image. These include \(\|\lmat{W}|\cdot|\|_p\) for certain matrices \(\lmat{W}\), Total Variation\cite{chambolle2004algorithm}, and multi-bang regularisation\cite{holman2020emission, richardson2021multi}.  The results for \(\|\lmat{W}|\cdot|\|_1\) and (isotropic) TV naturally coincide with those previously proven by G{\"u}ven \textit{et al}\cite{guven2016augmented}, though we also show that the former does not hold for any given matrix \(\lmat{W}\).  In section~\ref{sec: level set}, we also use our main result to provide a framework in which level set reconstruction methods\cite{osher2001level} may be applied directly to the magnitude image, allowing phase to vary arbitrarily spatially pixel-to-pixel.  To our knowledge level set methods have not previously been applied in this way, though in our primary consideration of SAR they have been applied for the segmentation of detected (i.e. magnitude only) imagery\cite{marques2011sar}.  We specifically apply our main theorem to the PaLEnTIR level set method\cite{ozsar2024parametric}, though it is broadly applicable, and may provide a basis to develop more efficient coherent level set reconstruction methods.

Secondly, in section~\ref{sec: gen funcs}, we provide a general method to numerically calculate \(\prox_{H(|\cdot|)}\) for such \(H\) which it cannot be guaranteed that our main result holds. We prove that, in this case, one needs to solve a proximal map of \(H\) bounded in the positive orthant, before applying the same phase correction as used in the first result.  A simple algorithm based on Douglas-Rachford splitting is proposed, providing a black-box approach to calculating \(\prox_{H(|\cdot|)}\) for any reasonable choice of regularisation function \(H\) (i.e. we can calculate \(\prox_H\) itself, and \(H\) is defined and somewhere finite on \(\mathds{R}^n_{\geq0}\)).  Since the starting iteration of this black-box algorithm applies the simple calculation of the first main result, wherever our first main result holds \(\prox_{H(|\cdot|)}\) is found without unnecessary Douglas-Rachford iterations.  Thus, this black box approach may always be applied to any such suitable regularisation function \(H\) without the need to determine if it meets sufficient conditions or otherwise, and without introducing unnecessary computational cost.

Finally, in section~\ref{sec: numerical results}, we apply these results to reconstructions of publicly available real airborne and satellite SAR data.  Applying our first result directly, we apply generalised Tikhonov and level-set reconstruction to the Gotcha carpark dataset\cite{casteel2007challenge}, including a multi-channel problem.  We then use our black-box general method to reconstruct data available from the Umbra open data program\cite{umbra} using Total Generalised Variation (TGV)\cite{bredies2010total}.
TGV can be shown not to satisfy the sufficient conditions of our primary result, therefore requiring this additional routine.  For SAR, TGV has only previously been applied to despeckling already formed magnitude imagery\cite{feng2015synthetic}, not as part of a reconstruction.

\section{Main result}\label{sec: main result}
Here we state the main result, which provides an efficient means to calculate \(\prox_{H(|\cdot|)}\) for certain functions \(H\), where \(|\cdot|:\mathds{C}^n\rightarrow\mathds{R}^n_{\geq0}\) is understood to mean element-wise absolute value of a vector throughout in a slight abuse of notation.

\begin{theorem}
	Let \(G(\lvec{z}) = H(|\lvec{z}|)\), \(G:\mathds{C}^n\rightarrow\mathds{R}\), where \(H:\mathds{R}^n\rightarrow \mathds{R}\) is a closed proper convex function. If \(\prox_H:\mathds{R}^n_{\geq0}\rightarrow \mathds{R}^n_{\geq0}\), then
	\begin{equation}
		\prox_G(\lvec{z}) \equiv \prox_H(\mathbf{r})\circ\lvec{\Phi},
		\label{eq: prox abs claim}
	\end{equation}
	where \(\lvec{r}:=|\lvec{z}|\), and \(\lvec{\Phi}:=\exp(\imi\angle\lvec{z})\), so that \(\lvec{z}\equiv\lvec{r}\circ\lvec{\Phi}\).
	\label{thm: main result}
\end{theorem}

Here, \(|\cdot|\) is understood to mean the element-wise absolute value of a vector, \(f(\lvec{z})=|\lvec{z}|:\mathds{C}^n\rightarrow\mathds{R}^{n}_{\geq0}\), and \(\circ\) denotes the element-wise (Hadamard) product.

\begin{proof}
	The claim (\ref{eq: prox abs claim}) can be written equivalently as
	\begin{equation}
		\prox_G(\lvec{z})\circ\overline{\lvec{\Phi}} \equiv \prox_H(\mathbf{r}),
	\end{equation}
	where \(\overline{\cdot}\) denotes complex conjugation, \(\overline{\lvec{\Phi}}\equiv \exp(-i\angle\lvec{z})\).  Expanding the left-hand side with the definitions of the proximal operator and \(G\), we have that 
	\begin{align}
		\prox_{G}(\lvec{z})\circ\overline{\lvec{\Phi}} &\equiv \left\{\argmin_{\lvec{y}\in\mathds{C}^n} G(\lvec{y}) + \frac{1}{2}\|\lvec{y}-\lvec{z}\|_3^2\right\}\circ\overline{\lvec{\Phi}} \\
		& \equiv \left\{\argmin_{\lvec{y}\in\mathds{C}^n} H(|\lvec{y}|) + \frac{1}{2}\|\lvec{y}-\lvec{z}\|_2^2\right\}\circ\overline{\lvec{\Phi}}.
	\end{align}
	Applying the change of variables \(\lvec{y}=\lvec{w}\circ\lvec{\Phi}\).  Then we have
	\begin{align}
		 \left\{\argmin_{\lvec{y}\in\mathds{C}^n} H(|\lvec{y}|) + \frac{1}{2}\|\lvec{y}-\lvec{z}\|_2^2\right\}\circ\overline{\lvec{\Phi}} &=  \argmin_{\lvec{w}\in\mathds{C}^n} H(|\lvec{w}\circ\lvec{\Phi}|) + \frac{1}{2}\|\lvec{w}\circ\lvec{\Phi}-\lvec{z}\|_2^2, \\
		 &= \argmin_{\lvec{w}\in\mathds{C}^n} H(|\lvec{w}|) + \frac{1}{2}\|\lvec{w}\circ\lvec{\Phi}-\lvec{z}\|_2^2,
		 \label{eq: almost proxH}
	\end{align}
since applying the changes of phase \(\lvec{w}\circ\lvec{\Phi}\) does not affect its absolute values. We can also write
\begin{align}
	\|\lvec{w}\circ\lvec{\Phi}-\lvec{z}\| =& \|\overline{\lvec{\Phi}}\|\|\lvec{w}\circ\lvec{\Phi}-\lvec{z}\|\nonumber\\
	=&  \|\overline{\lvec{\Phi}}\circ(\lvec{w}\circ\lvec{\Phi}-\lvec{z})\| \nonumber\\
	=& \|\lvec{w}-\lvec{z}\circ\overline{\lvec{\Phi}}\|= \|\lvec{w}-\lvec{r}\|.	
\end{align}
Since \(\lvec{r}\in\mathds{R}^n_{\geq0}\subset \mathds{R}^n\), we have
\begin{equation}
	\|\lvec{w}-\lvec{r}\|\geq \|\lvec{x}-\lvec{r}\|, \quad \lvec{x}=\Re(\lvec{w}), \quad \forall \lvec{w}\in\mathds{C}^n.
\end{equation}
It is now straightforward to see that
\begin{align}
	\argmin_{\lvec{w}\in\mathds{C}^n}H(|\lvec{w}|) + \frac{1}{2}\|\lvec{w}\circ\lvec{\Phi} - \lvec{z}\|^2_2 =& \argmin_{\lvec{x}\in\mathds{R}^n}H(|\lvec{x}|) + \frac{1}{2}\|\lvec{x}-\lvec{r}\|_2^2, \\
	\equiv&\argmin_{\lvec{x}\in\mathds{R}^n_{\geq0}}H(\lvec{x}) + \frac{1}{2}\|\lvec{x}-\lvec{r}\|_2^2.
	\label{eq: equiv +ve orth}
\end{align}
Since by premise \(\prox_H(\lvec{r})\in\mathds{R}^n_{\geq0}\) \(\forall\lvec{r}\in\mathds{R}^n_{\geq0}\), and by definition we have \(\lvec{r}\in\mathds{R}^n_{\geq0}\),
\begin{equation}
	\prox_H(\lvec{r}):=\argmin_{\lvec{x}\in\mathds{R}^n} H(\lvec{x})+\frac{1}{2}\|\lvec{x}-\lvec{r}\|_2^2\equiv \argmin_{\lvec{x}\in\mathds{R}_{\geq 0}^n}H(\lvec{x})+\frac{1}{2}\|\lvec{x}-\lvec{r}\|_2^2
	\label{eq: equiv final}
\end{equation}
and we have the required result.
\end{proof}
This sufficient appears to be as general as possible, or the \textit{weakest} sufficient condition in the sense of Lin\cite{lin2001strongest}. That is, we can only guarantee the final equivalence in (\ref{eq: equiv final}) if we guarantee \(\prox_H(\lvec{r})\in\mathds{R}^n_{\geq0}\) \(\forall\lvec{r}\in\mathds{R}^n_{\geq0}\).  This of course may also happen to occur where it is not guaranteed for all such \(\lvec{r}\).  Indeed, we can see that the only reason this is not also a necessary condition and the proof does not run in reverse is that we do not require \(\prox_H(\lvec{r})\in\mathds{R}^n_{\geq}\forall\mathds{r}\in\mathds{R}_{\geq0}^+\), but only pointwise for the specific \(\mathds{r}\) at which we evaluate the map. Similarly, \(H\) is assumed to be closed, proper and convex to guarantee a unique solution to the proximal map\cite{parikh2014proximal} -- though one may exist in cases where this does not hold.

It leads naturally to a simple corollary.

\begin{cor}
	Let \(G(\lvec{z}) = H(|\lvec{z}|)\), \(G:\mathds{C}^n\rightarrow\mathds{R}\), where \(H\) is a closed proper convex function with effective domain \(\dom H \subseteq \mathds{R}_{\geq0}^n\).  Then
	\begin{equation}
		\prox_G(\lvec{z}) \equiv \prox_H(\mathbf{r})\circ\lvec{\Phi},
	\end{equation}
	where \(\lvec{r}:=|\lvec{z}|\), and \(\lvec{\Phi}:=\exp(\imi\angle\lvec{z})\).
	\label{cor: eff dom} 
\end{cor}
Clearly \(\dom H\subseteq\mathds{R}^n_{\geq0}\Rightarrow \prox_H(\lvec{r})\in\mathds{R}_{\geq0}^n\ \forall \lvec{r}\). The proof also follows directly from (\ref{eq: equiv +ve orth}), since extending the minimisation to be outside of the effective domain of \(H\) cannot affect the solution.

A second corollary may be more practically applicable.
\begin{cor}
	Let \(G(\lvec{z}) = H(|\lvec{z}|)\), \(G:\mathds{C}^n\rightarrow\mathds{R}\), where \(H:\mathds{R}^n\rightarrow \mathds{R}\) is a closed proper convex function with the property \(H(|\lvec{x}|)\leq H(\lvec{x})\ \forall \lvec{x}\in\mathds{R}^n\). Then
	\begin{equation}
		\prox_G(\lvec{z}) \equiv \prox_H(\mathbf{r})\circ\lvec{\Phi},
		\label{eq: prox abs claim 2}
	\end{equation}
	where \(\lvec{r}:=|\lvec{z}|\), and \(\Phi:=\exp(\imi\angle\lvec{z})\), so that \(\lvec{z}\equiv\lvec{r}\circ\Phi\).
	\label{cor: H bounds}
\end{cor}
\begin{proof}
We begin once again at (\ref{eq: equiv +ve orth}). To extend the minimisation to the whole of \(\mathds{R}^n\) for the required result, we need that
\begin{align}
	H(\lvec{x}) +\frac{1}{2}\|\lvec{x}-\lvec{r}\|_2^2\geq& H(|\lvec{x}|) +\frac{1}{2}\||\lvec{x}|-\lvec{r}\|_2^2 \quad \forall \lvec{x}\in\mathds{R}^n, \nonumber \\
	\Leftrightarrow \quad 	H(\lvec{x})\geq& H(|\lvec{x}|) +\frac{1}{2}\||\lvec{x}|-\lvec{r}\|_2^2 - \frac{1}{2}\|\lvec{x}-\lvec{r}\|_2^2.
\end{align}
To see this, 
\begin{align}
	H(|\lvec{x}|) + \frac{1}{2}\||\lvec{x}|-\lvec{r}\|_2^2 - \frac{1}{2}\|\lvec{x}-\lvec{r}\|_2^2
	=& H(|\lvec{x}|) + \frac{1}{2}\left((|\lvec{x}|-\lvec{r}-\lvec{x}+\lvec{r})\cdot(|\lvec{x}|-\lvec{r}+\lvec{x}-\lvec{r})\right)\nonumber\\
	=& H(|\lvec{x}|) -\lvec{r}\cdot(|\lvec{x}|-\lvec{x}) \leq H(|\lvec{x}|)
	\label{eq: H bounds}
\end{align}
since \(\lvec{r}\in\mathds{R}_{\geq0}^n\). By proposition \(H(\lvec{x})\geq H(|\lvec{x}|)\), so the result holds.
\end{proof}
~\\
\begin{remark}
	It is clear that the result (\ref{eq: prox abs claim}) will also hold for any \(H\) and \(\lvec{r}=|\lvec{z}|\) such that
	\begin{equation}
		H(\lvec{x}) \geq H(|\lvec{x}|) - \lvec{r}\cdot(|\lvec{x}|-\lvec{x}) \quad \forall \lvec{x}\in\mathds{R}^n,
	\end{equation}
	which follows directly from (\ref{eq: H bounds}). However, algorithmically this may not be straightforward (or useful) to enforce.
\end{remark}~\\

\section{Application to some common regularisation functions}\label{sec: common regs}
Here we show how the main results of Section~\ref{sec: main result} may be applied to guarantee a simple proximal map to certain common regularisation functions applied directly to the magnitude image.  In each case, we provide only a sketch proof that they meet the conditions above.

\subsection{Vector norms}
From Corollary~\ref{cor: H bounds}, we are able to deduce the proximal map of several functions commonly used as regularisation terms when applied to the magnitude of a complex image.  The first, somewhat trivially, is weighted p-norms
\begin{equation}
	\prox_{\||\cdot|\|_{\lmat{W},p}}(\lvec{z}) = \prox_{\|\cdot \|_{\lmat{W},p}}(\lvec{r})\circ\lvec{\Phi}\quad \forall p\geq1, \quad\lmat{W}=\diag(\lvec{w})
\end{equation}
where \(\lvec{r}=|\lvec{z}|\) and \(\lvec{\Phi}=\exp{\imi\angle\lvec{z}}\).  This follows from Corollary~\ref{cor: H bounds} since flipping signs of vector elements does not effect their \(L_p\) norm. 

In several cases we can include a non-diagonal matrix \(\lmat{W}\) in the p-norm,
\begin{equation}
	\prox_{\|\lmat{W}|\cdot|\|_p}(\lvec{z}) = \prox_{\|\lmat{W}\cdot\|_p}(\lvec{r}).
	\label{eq: prox WLp}
\end{equation}
We can show that (\ref{eq: prox WLp}) holds for any \(\lmat{W}\in\mathds{R}^{N\times N}\) and \(p\geq1\) such that either
\begin{align}
	\lmat{W}=& \bigl[w_1(\lvec{e}_{i_{11}}-\lvec{e}_{i_{12}}), \ldots, w_N(\lvec{e}_{i_{N1}}-\lvec{e}_{i_{N2}})\bigr]^T,\mbox{ or} \label{eq: difference W}\\
	\lmat{W}^T\lmat{W}=&\lmat{I}. \label{eq: orthogonal W}
\end{align}
For a weighted difference matrix (\ref{eq: difference W}) we must have \(\|\lmat{W}|\lvec{x}|\|_p \leq \|\lmat{W}\lvec{x}\|_p\) by the reverse triangle inequality, i.e.
\begin{equation}
	\bigl|\|u\|-\|v\|\bigr|\leq \|u-v\|,
	\label{eq: reverse triangle}
\end{equation} 
so can again apply Corollary~\ref{cor: H bounds}.

In the case of an orthogonal matrix (\ref{eq: orthogonal W}) we can apply the equivalence\cite{parikh2014proximal}
\begin{equation}
	\prox_{\|\lmat{W}\cdot\|_p}(\lvec{v})\equiv \lmat{W}^T\prox_{\|\cdot\|_p}(\lmat{W}\lvec{v}).
\end{equation}
We can see that \(\lvec{x}:=\prox_{\|\cdot\|_p}(\lmat{W}\lvec{v})\) must lie in the same orthant as \(\lmat{W}\lvec{v}\), and since \(\lmat{W}\cdot\) is an orthogonal transformation this is the orthant described by the \(\lmat{W}\lvec{e}_i\) directions.  Then, again using the angle-preserving property of orthogonal matrices, \(\lmat{W}^T\lvec{x}\) must lie in the orthant in the \(\lmat{W}^T\lmat{W}\lvec{e}_i\) directions, i.e. the positive orthant.  Therefore, for an orthogonal matrix, the map \(\prox_{\|W\cdot\|_p}(\lvec{v})\in\mathds{R}^n_{\geq0}\) \(\forall \lvec{v}\in\mathds{R}^n_{\geq0}\), and we can apply the result of Theorem~\ref{thm: main result} directly.

\begin{remark}
	An alternative proof for (\ref{eq: prox WLp}) has previously been provided by G{\"u}ven \textit{et al}\cite{guven2016augmented} in the case of a 1-norm, but without restriction on \(\lmat{W}\).  However, we can see that this does not hold for any \(\lmat{W}\) by providing a simple counterexample, for which it's sufficient to show a proximal map for real-valued \(\lvec{z}\).  Taking \(\lmat{W}=\left[\begin{smallmatrix}
		1.  & -0.7  &  0.35\\
		-0.7 &   1. &  -0.9\\
		0.35 &  -0.9  &  1.
	\end{smallmatrix}\right]\), and \(\lvec{z} = [2,\ 10^{-9},\ 10^{-9}]^T\), solving numerically we find that (to three significant figures),
	\begin{align}
		\prox_{\|\lmat{W}|\cdot|\|_1}(\lvec{z}) = & [0.815, 0.576, 0.005]^T,\\
		\prox_{\|\lmat{W}\cdot\|_1}(|\lvec{z}|)\circ\exp(\imi\angle\lvec{z})\equiv \prox_{\|\lmat{W}\cdot\|_1}(|\lvec{z}|) = & [0.826, 0.555, -0.025]^T.
		\label{eq: counter example}
	\end{align}
	So we see it is possible not just for the result to be out of phase (in this case an element \qty{\pi}{rad} out of phase), but also for the result to have the incorrect magnitudes of elements.  
\end{remark}

\subsection{Total Variation and variants}
A common choice of regularisation term is the Total Variation (TV) semi-norm\cite{chambolle2004algorithm, beck2009ieee}, which promotes piecewise-constant solutions to the inverse problem (i.e. sparse in gradient), defined as isotropic and anisotropic TV for a 2D image as either
\begin{align}
	\TV_i(\lvec{u}):=& \|\lmat{D} \lvec{u}\|_{2,1} = \sum_{i,j} \sqrt{|(\lmat{D}_x \lvec{\lvec{u}})_{i,j}|^2 + |(\lmat{D}_y \lvec{\lvec{u}})_{i,j}|^2} \qquad &\mbox{(isotropic)},
	\label{eq: iso TV}\\
	\TV_a(\lvec{u}):=& \|\lmat{D} \lvec{u}\|_{1} = \sum_{i,j}\bigl\{ |(\lmat{D}_x \lvec{\lvec{u}})_{i,j}| + |(\lmat{D}_y \lvec{\lvec{u}})_{i,j}|\bigr\} \qquad &\mbox{(anisotropic)},
	\label{eq: aniso TV}
\end{align}
respectively, or equivalently in 3D as
\begin{align}
	\TV_i(\lvec{u}):=& \|\lmat{D} \lvec{u}\|_{2,1} = \sum_{i,j,k}\sqrt{|(\lmat{D}_x \lvec{\lvec{u}})_{i,j,k}|^2 + |(\lmat{D}_y \lvec{\lvec{u}})_{i,j,k}|^2+|(\lmat{D}_z\lvec{\lvec{u}})_{i,j,k}|^2},
	\label{eq: iso TV}\\
	\TV_a(\lvec{u}):=& \|\lmat{D} \lvec{u}\|_{1} = \sum_{i,j,k}\bigl\{ |(\lmat{D}_x \lvec{\lvec{u}})_{i,j,k}| + |(\lmat{D}_y \lvec{\lvec{u}})_{i,j,k}|+|(\lmat{D}_z \lvec{\lvec{u}})_{i,j,k}|\bigr\},
	\label{eq: aniso TV}
\end{align}
where \(\lmat{D}\) is the discrete gradient operator, and \(\lmat{D}_x\) \(\lmat{D}_y\) and \(\lmat{D}_z\) are the discrete first derivative operators in the \(x\), \(y\) and \(z\) image coordinate directions, with the subscript \(i,j,k\) in referring to pixels/voxels these image coordinate directions. For multi-channel or hyperspectral reconstruction problems, one might also use vectorial (channelwise) TV (VTV), defined in 2D for example as \cite{duran2016collaborative, papoutsellis2021core2}
\begin{equation}
	\VTV(\lvec{u}):=\|D\lvec{u}\|_{2,1} = \sum_{k}\sum_{i,j}\sqrt{|(\lmat{D}_x \lvec{\lvec{u}}_k)_{i,j}|^2 + |(\lmat{D}_y \lvec{\lvec{u}}_k)_{i,j}|^2 },
	\label{eq: VTV}
\end{equation}
where \(\lvec{u}:=[\lvec{u}_1,\ldots,\lvec{u}_K]^T\) is the multichannel image.  In other cases one may wish to include the derivative across the channels or hyperspectral component. For example, where these represent sequential timesteps of data, one may use spatio-temporal TV\cite{papoutsellis2021core2, watson2022focusing}
\begin{equation}
	\TV_{st}(\lvec{u}):= \|\lmat{D} \lvec{u}\|_{2,1} = \sum_{t,i,j}\sqrt{|(\lmat{D}_t\lvec{\lvec{u}})_{t,i,j}|^2 + |(\lmat{D}_x \lvec{\lvec{u}})_{t,i,j}|^2 + |(\lmat{D}_y \lvec{\lvec{u}})_{t,i,j}|^2}.
	\label{eq: spatiotemporal TV}
\end{equation}

As before mentioned, G\"{u}ven et al have previously provided a technical proof that \(\prox_{\TV}\) follows the equivalence (\ref{eq: prox abs claim}) for isotropic TV in 2D\cite{guven2016augmented}.  We can also succinctly show that for both the iso- and anisotropic TV in either 2 or 3 dimensions, as well as the multi-channel cases, all meet the requirements of Corollary~\ref{cor: H bounds}: once again using the reverse triangle inequality, we must have \(\TV(|\lvec{x}|)\leq \TV(\lvec{x})\) \(\forall \lvec{x}\in\mathds{R}^n\).  Thus, the result of Theorem~\ref{thm: main result} holds.

\subsection{Restriction to feasible sets}
Often one wishes to restrict the solution to some feasible set of values, \(\lvec{x}\in\mathcal{C}\).  For example this may be to set the solution to zero outside of a region (e.g. outside of a body in medical imaging), or for the values to lie within some feasible bounds.  Carrying out a bounded reconstruction problem over \(\mathcal{C}\) is equivalent to adding the indicator function \(\chi_{\mathcal{C}}\) as a regularisation term,
\begin{equation}
	\chi_{\mathcal{C}}(\lvec{x}):=\left\{\begin{array}{ll}
		0 & \mbox{for } \lvec{x}\in\mathcal{C},\\
		\infty& \mbox{otherwise}.
		\end{array}\right.
\end{equation}
This has proximal map given by the projection onto \(\mathcal{C}\),
\begin{equation}
	\prox_{\chi_{\mathcal{C}}}(\lvec{x})\equiv \proj_{\chi_{\mathcal{C}}}(\lvec{x}).
\end{equation}
Clearly if \(\mathcal{C}\subseteq\mathds{R}^n_{\geq}\) then we may apply Corollary~\ref{cor: eff dom}.

One may also wish to constrain the reconstruction to taking only a finite set of values \(\mathcal{A}=\{a_1,\ldots a_k\}\), referred to as being \textit{multi-bang}\cite{holman2020emission, richardson2021multi}. This can be enforced by the regularisation penalty term
\begin{align}
	\mathcal{M}(\lvec{x}):=&\sum_i m(\lele{x}_i), \\
	m(x)=&\left\{\begin{array}{ll}(a_{i+1}-x)(x-a_i),\quad& \mbox{if }x\in[a_i,a_{i+1}]\mbox{ for some } i,\\
	\infty,& \mbox{otherwise},\end{array}\right.
\end{align}
where the \(a_i\) are ordered in increasing value.  This non-convex penalty term has proximal map which has multi-bang values as stationary points given by\cite{holman2020emission}
\begin{equation}
	[\prox_{\frac{1}{\tau}\mathcal{M}}(\lvec{x})]_j = \left\{\begin{array}{ll} a_0  &   \mbox{if } x_j\leq x_{0,+},\\
		a_i & \mbox{if } x_{i,-}\leq x_j \leq x_{i,+}\mbox{ for } i\in\{1,\ldots,{k-1}\},\\
		a_n & \mbox{if } x_{n,-}\leq x_j, \\
		\frac{1}{1-2\tau}(x) & \mbox{if } x_{i,+}<x_j<x_{i+1,-}\mbox{ for } i\in\{0,\ldots,{k-1}\}
	\end{array}\right.
\end{equation}
where
\begin{align*}
	x_{i,-}=&a_i - \tau(a_i - a_{i-1})\mbox{ for }i=1,\ldots,k\\
	x_{i,+}=& a_i+\tau(a_{i+1}-a_i)\mbox{ for }i=0,\ldots,n-1.
\end{align*}
Clearly, if the target values \(a_i\geq0\) then \(\prox_{\frac{1}{\tau}\mathcal{M}}(\lvec{x})\in\mathds{R}^n_{\geq0}\) \(\forall \mathds{x}\in\mathds{R}^n_{\geq0}\) for admissible \(0<\tau<\tfrac{1}{2}\), so for a suitable choice of regularisation parameter (and step length in a chosen optimisation algorithm), despite not being a convex function, the result of Theorem~\ref{thm: main result} holds.

Multi-bang may not be directly applicable to our main motivation of SAR, though may be useful in other coherent imaging problems such as electromagnetic imaging in medical settings.

\section{Application to Level Set Reconstruction}\label{sec: level set}
In the level set approach, one formulates the problem in the form
\begin{equation}
	\min_{\lvec{p}} \frac{1}{2}\|\mathcal{M}(f(\lvec{r};\lvec{p})) - \lvec{d}\|_2^2,
	\label{eq: level set}
\end{equation}
where \(\mathcal{M}\) is the forward operator, and \(f\) maps parameters describing the level set to the (discrete) values in the image.    The level set reconstruction problem (\ref{eq: level set}) can be equivalently written as
\begin{equation}
\min_{\lvec{x}\in \mathcal{C}} \|\mathcal{M}(\lvec{x}) - \lvec{d}\|_2^2 \equiv \min_{\lvec{x}\in\mathds{R}^n} \|\mathcal{M}(\lvec{x}) - \lvec{d}\|_2^2 + \chi_{\mathcal{C}}(\lvec{x}),
\label{eq: equiv level set}
\end{equation}
where \(\mathcal{C}\subset\mathds{R}^n\) is the space of functions which can be represented by the level set basis functions and \(\chi_{\mathcal{C}}\) is the indicator function for the subspace \(\mathcal{C}\subset\mathds{R}^n\). 

If \(\mathcal{C}\subseteq\mathds{R}^n\), then \(\dom(f)\subseteq\mathds{R}^n\) so clearly (\ref{eq: equiv level set}) meets the requirements of Corollary~\ref{cor: eff dom}.  For example, consider the PaLEnTIR formulation of Ozsar \textit{et al}\cite{ozsar2024parametric} based on radial basis functions \(\phi\),
\begin{align}
	f(\lvec{r};\lvec{p}, C_H, C_L) =& C_H T_w(\phi(\lvec{r};\lvec{p}) - c) + C_L(1-T_w(\phi(\lvec{r};\lvec{p})-c)) \in \mathds{R}^n, \label{eq: palentir f} \\
	\phi(\lvec{r};\lvec{p}) =& \sum_{j=1}^N \sigma_h(\alpha_i)\psi(\lmat{R}_j(\lvec{r}-\lvec{\chi}_j))\quad \lvec{p}^T:=[\lvec{\alpha}^T, \lvec{\beta}^T, \lvec{\gamma}^T],
\end{align}
with matrices \(\lmat{R}_j:=\lmat{R}(\beta_j,\gamma_j)\) depending on subvectors of \(\lvec{\beta}\) and \(\lvec{\gamma}\) respectively, \(T_w\) is a transition function giving a smooth approximation to the Heaviside step function, and where \(c\) defines the level set to be taken (typically \(0\leq c\ll 1\)).  Thus \(\dom(f)\) is given by the set
\begin{equation}
	\begin{aligned}
	\mathcal{C}_{pal}:=&\left\{\lvec{x}\,|\, \lvec{x}=C_H T_w(\phi(\lvec{r};\lvec{p}) - c) + C_L(1-T_w(\phi(\lvec{r};\lvec{p})-c))  \right\} \equiv \Span\{f(\lvec{p})\}
	\end{aligned}
	\label{eq: palentir set C}
\end{equation}
This provides us with the following Lemma.

\begin{lemma}
	For any \(C_H> C_L\geq 0\), defining \(G(\lvec{z}):=\chi_{\mathcal{C}_{pal}}(|\lvec{z}|)\), then
	\begin{equation}
		\prox_G(\lvec{z})\equiv \prox_{\chi_{\mathcal{C}_p}}(\lvec{r})\circ\lvec{\Phi},
		\label{eq: prox level set}
	\end{equation}
	where \(\lvec{r}=|\lvec{z}|\) and \(\lvec{\Phi} = \exp{i\angle\lvec{z}}\). Moreover, 
	\begin{equation}
		\prox_{\chi_{\mathcal{C}_{pal}}}(\lvec{r}) = f(\tilde{\lvec{p}}), \qquad \tilde{\lvec{p}} = \argmin_{\lvec{p}} \|f(\lvec{p}) - \lvec{r}\|_2^2.
		\label{eq: level set image denoise}
	\end{equation}
\label{lem: prox palentir}
\end{lemma}
\begin{proof}
	From definition (\ref{eq: palentir set C}) it is clear that \(f\geq0\) if \(C_H,C_L\geq0\), so \(\mathcal{C}_p\subseteq\mathds{R}^n_{\geq0}\), so (\ref{eq: prox level set}) is obtained from Lemma~\ref{cor: eff dom}.  The proximal map (\ref{eq: level set image denoise}) simply arises from definition,
	\begin{align}
		\prox_{\chi_{\mathcal{C}_{pal}}}(\lvec{r})=&\proj_{\chi_{\mathcal{C}_{pal}}}(\lvec{r}) \nonumber\\
		=& \argmin_{\lvec{x}\in\mathcal{C}_{pal}}\frac{1}{2}\|\lvec{x}-\lvec{r}\|_2^2,
	\end{align}
\begin{equation}
	\Leftrightarrow\quad \lvec{x}=f(\tilde{\lvec{p}}),\qquad	\tilde{\lvec{p}}= \argmin_{\lvec{p}}\frac{1}{2}\|f(\lvec{p})-\lvec{r}\|_2^2.
\end{equation}
\end{proof}

An equivalent result to Lemma~\ref{lem: prox palentir} could be simply created for any level set formulation using suitable basis functions and parameterisation which are positive valued.  Use of (\ref{eq: prox level set}) and (\ref{eq: level set image denoise}) in a reconstruction procedure to solve (\ref{eq: equiv level set}) reconstruction method reduces the level set reconstruction problem to a simpler image de-noising one which generally should be inexpensive to solve.  A potential drawback that the level set method does not ``know'' anything about the forward operator \(\mathcal{M}\) (i.e. no information about derivatives of \(\mathcal{M}\), so a simple proximal gradient algorithm may require many iterations.  Nevertheless, this provides a mechanism for application of level set reconstruction techniques to the magnitude only of complex-valued images, allowing phase to vary pixel-by-pixel within each constant-valued region.

\section{Proximal map for more general functions}\label{sec: gen funcs}
For many functions the sufficient condition of Theorem~\ref{thm: main result} may not hold, and as seen in the example (\ref{eq: counter example}) this may result in the equality (\ref{eq: prox abs claim}) not holding.  In such cases, we need to solve a bounded proximal map for \(H\), for which we use the notation \(\prox^+\)
\begin{equation}
	\begin{aligned}
	\prox_H^+(\lvec{r}):=&\argmin_{\lvec{x}\in\mathds{R}^n_{\geq0}}H(\lvec{x}) + \frac{1}{2}\|\lvec{x}-\lvec{r}\|_2^2 \\
	\equiv& \argmin_{\lvec{x}}H(\lvec{x}) + \frac{1}{2}\|\lvec{x}-\lvec{r}\|_2^2 + \chi_{\mathds{R}^n_{\geq0}} \\
	\equiv& \prox_{H+ \chi_{\mathds{R}^n_{\geq0}}}.
	\end{aligned}
\label{eq: bounded prox}
\end{equation}
Since \(\mathds{R}^n_{\geq0}\) is a convex set and \(H(\lvec{x})+\tfrac{1}{2}\|\lvec{x} - \lvec{r}\|_2^2\) is (strictly) convex, the solution to (\ref{eq: bounded prox}) exists and is stable.  Thus we are able to solve (\ref{eq: bounded prox}) numerically (assuming a closed form solution does not exist) via a suitable optimisation scheme such as Douglas-Rachford splitting, PDHG or ADMM.

Therefore, if we're able to solve (\ref{eq: bounded prox}) efficiently, then we are able to calculate the proximal map of a function \(H(|\lvec{z}|)\) via the following Theorem.

\begin{theorem}
	Let \(G(\lvec{z})=H(|\lvec{z}|)\), \(G:\mathds{C}^n\rightarrow\mathds{R}\), where \(\lvec{H}:\mathds{R}^n\rightarrow\mathds{R}\) is a closed proper convex function. Then
	\begin{equation}
		\prox_G(\lvec{z})\equiv\prox_H^+(\lvec{r})\circ\lvec{\Phi},
	\end{equation}
where \(\lvec{r}:=|\lvec{z}|\) and \(\lvec{\Phi}:=\exp(i\angle\lvec{z})\), so that \(\lvec{z}\equiv \lvec{r}\circ\lvec{\Phi}\).
\label{thm: bounded equiv}
\end{theorem}
The proof of Theorem~\ref{thm: bounded equiv} is the same as for Theorem~\ref{thm: main result} up to equation (\ref{eq: equiv +ve orth}) which is the required result.  Note that in this case, Theorem~\ref{thm: bounded equiv} provides both a necessary and sufficient condition.

Algorithm~\ref{alg: prox abs} provides a simple scheme to solve (\ref{eq: bounded prox}) via Douglas-Rachford splitting to minimise \(H(\lvec{x})+F(\lvec{x})\), given by
\begin{subequations}
	\begin{align}
		\lvec{x}_{k+1} =& \prox_H(\lvec{y}_k) \\
		\lvec{y}_{k+1} =& \lvec{y}_k + \prox_F(2\lvec{x}_{k+1}-\lvec{y}) - \lvec{x}_{k+1},
	\end{align}
\end{subequations}
 by taking \(F(\lvec{x}):=\chi_{\mathds{R}^n_{\geq0}}+\frac{1}{2}\|\lvec{x}-\lvec{r}\|_2^2\).  This requires the proximal map of \(F(\lvec{x})\), given by
\begin{align}
	\prox_{F}(\lvec{y}):=&\argmin_{\lvec{x}} \chi_{\mathds{R}^n_{\geq0}} + \frac{1}{2}\|\lvec{x}-\lvec{r}\|_2^2 + \frac{1}{2}\|\lvec{x}-\lvec{y}\|_2^2\label{eq: prox chi plus}\\
	\equiv &\argmin_{\lvec{x}} \chi_{\mathds{R}^n_{\geq0}} + \frac{1}{2}\bigl\|\lvec{x} - \left(\frac{\lvec{r}+\lvec{y}}{2}\right)\bigr\|_2^2 \label{eq: intermed prox chi}  \\
	=& \proj_{\mathds{R}^n_{\geq0}}\left(\frac{\lvec{r}+\lvec{y}}{2}\right) = \prox_{\chi_{\mathds{R}^n_{\geq0}}}\left(\frac{\lvec{r}+\lvec{y}}{2}\right),
\end{align}  
with equivalence of (\ref{eq: intermed prox chi}) seen by expanding the least-squares term of (\ref{eq: prox chi plus}) and completing the square (which does not change the location of the minimum).  

\begin{algorithm}[h]
	\caption{Complex proximal-abs}\label{alg: prox abs}
	\begin{algorithmic}
		\State \(\lvec{r}\leftarrow|\lvec{z}|\)
		\State \(\lvec{\Phi}\leftarrow \exp(\imi\angle\lvec{z})\)
		\State \(\lvec{y}\leftarrow\lvec{r}\)
		\State \(\lvec{x}\leftarrow\prox_{H}(\lvec{y})\)
		\While{{\bf any} $\lele{x}_i<0$}
			\State \(\lvec{y}\leftarrow \lvec{y}+\proj_{\mathds{R}^n_{\geq0}}(\lvec{x}-0.5\lvec{y}+0.5\lvec{r}) - \lvec{x}\)
			\State \(\lvec{x}\leftarrow\prox_H(\lvec{y})\)
		\EndWhile
		\State \(\tilde{\lvec{z}}\leftarrow\lvec{x}\circ\lvec{\Phi}\)
		\State {\bf Return} \(\tilde{\lvec{z}}\) \Comment{\(\tilde{\lvec{z}}=\prox_{H(|\cdot|)}(\lvec{z})\)}
	\end{algorithmic}
\end{algorithm}

The initial choice of \(\lvec{y}_0\) in Algorithm~\ref{alg: prox abs} ensures that any function \(H\) satisfying the condition of Theorem~\ref{thm: main result} avoids unnecessary Douglas-Rachford iterations.  This simple base algorithm should also have a suitable guard against too many iterations. 

\FloatBarrier
\section{Application to SAR data}\label{sec: numerical results}
Here we apply the results of the previous sections to the regularisation of SAR imagery.  The purpose is to demonstrate that Algorithm~\ref{alg: prox abs} may readily be applied to a wide variety of real-data SAR imaging problems, providing a simple numerical recipe for different forms of regularisation -- in particular some examples which have not previously been applied to coherent SAR imaging. For the least-squares SAR reconstruction problem
\begin{equation}
	\minimise_{\lvec{x}\in\mathds{C}^n} \frac{1}{2}\|\lmat{A}\lvec{x}-\lvec{d}\|_2^2 + \lambda R(\lvec{x}),\quad\lambda>0,
	\label{eq: gen reconstruction}
\end{equation}
the cases of generalised Tikhonov regularisation is shown in section~\ref{subsec: gen tik}, level set reconstruction in section~\ref{subsec: level set}, and total generalised variation (TGV) in section~\ref{subsec: tgv} for regularisation function \(R\).  Details of the SAR data model, \(\lmat{A}\), are provided in Appendix~\ref{ap: data model}.

For these reconstructions, we use publicly available data from both the Gotcha challenge set\cite{casteel2007challenge} and Umbra open data program\cite{umbra}.  The reconstructions themselves are carried out using routines available in the CCPi Core Imaging Library (CIL)\cite{jorgensen2021core1,papoutsellis2021core2}, and specifically the implementation of the \texttt{PDHG} algorithm\cite{chambolle2011first} therein.  This solves optimisation problems of the form
\begin{equation}
	\min_{\lvec{x}}f(\lmat{K}\lvec{x}) + g(\lvec{x})
\end{equation}
where \(f\) and \(g\) are a convex functions with ``simple'' proximal maps.  This allows us to consider problems involving \(\tilde{f}(\lvec{x}):=f(\lmat{K}\lvec{x})\) where \(\tilde{f}\) itself does not have an easy to calculate proximal map.  For each of the reconstructions provided, we use the splitting \(f = 0.5\|\cdot - \lvec{d}\|_2^2\), \(\lmat{K}=\lmat{A}\) in PDHG. Regularisation is included in \(g\), since the nonlinearity of applying regularisation to the magnitude precludes a simple formulation in the form \(\tilde{f}\) beyond taking \(\lmat{K}\) as the identity.

We make available an implementation of Algorithm~\ref{alg: prox abs} which may be used with arbitrary regularisers (a \texttt{Function} in the CIL class structure).  This takes as input any other CIL \texttt{Function}, provided they are defined and bounded on (a subset of) \(\mathds{R}^n_{\geq0}\).

\subsection{Multi-look SAR with generalised Tikhonov}\label{subsec: gen tik}
Here, we apply generalised Tikhonov regularisation to a joint reconstruction multi-look SAR imagery.  That is, multiple separate spotlight SAR collections taken from different observation angles in azimuth. Letting \(\lvec{z}:=[\lvec{z}_1,\ldots,\lvec{z}_n]^T\) for \(\lvec{z}_i\) the \(n\) individual ``single-look'' image channels, then we solve
\begin{equation}
	\hat{\lvec{z}} = \argmin_{\lvec{z}} \sum_{i=1}^n\left\{\frac{1}{2}\|\lmat{A}_i\lvec{z}_i - \lvec{d}_i\|_2^2 \right\} + \frac{\lambda}{2}\|\lmat{D}|\lvec{z}|\|_2^2,
	\label{eq: gen tik}
\end{equation}
where \(\lmat{D}\) is discrete gradient in both space and image channel (ordered by look angle), \(\lmat{A}_i\) is the forward operator for image channel \(i\) and \(\lvec{d}_i\) the respective measured data, and \(\lvec{z}:=[\lvec{z}_1,\ldots,\lvec{z}_n]\) the combined set of reconstructed images for each look direction/channel.   The proximal of (\ref{eq: gen tik}), \(\prox_{\frac{\lambda}{2}\|\lmat{D}\cdot\|_2^2}(\lvec{r})\) is itself calculated numerically with 50 iterations of the PDHG algorithm, using the splitting \(f=\frac{1}{2}\|\cdot\|_2^2\), \(g=\frac{1}{2\lambda}\|\lvec{r}-\cdot\|_2^2\), and \(\lmat{K}=\lmat{D}\) the discrete gradient in both space and look/channel. Since there is not a natural scaling between this channel and pixel spacing, we set the channel-wise discrete derivative to 10 times the spatial one for the purpose of an example.  

This should promote magnitude images which vary smoothly with channel, having a similar effect to incoherent combination of the subaperture images, for example as in Stevens \textit{et al}\cite{stevens2017bright}, whilst retaining stronger variations between each channel and varying smoothly spatially.

\begin{wrapfigure}{R}{0.33\textwidth}
	\centering
	\scalebox{-1}[1]{\includegraphics[width=0.28\textwidth, trim=1cm 1cm 1cm 1cm, clip]{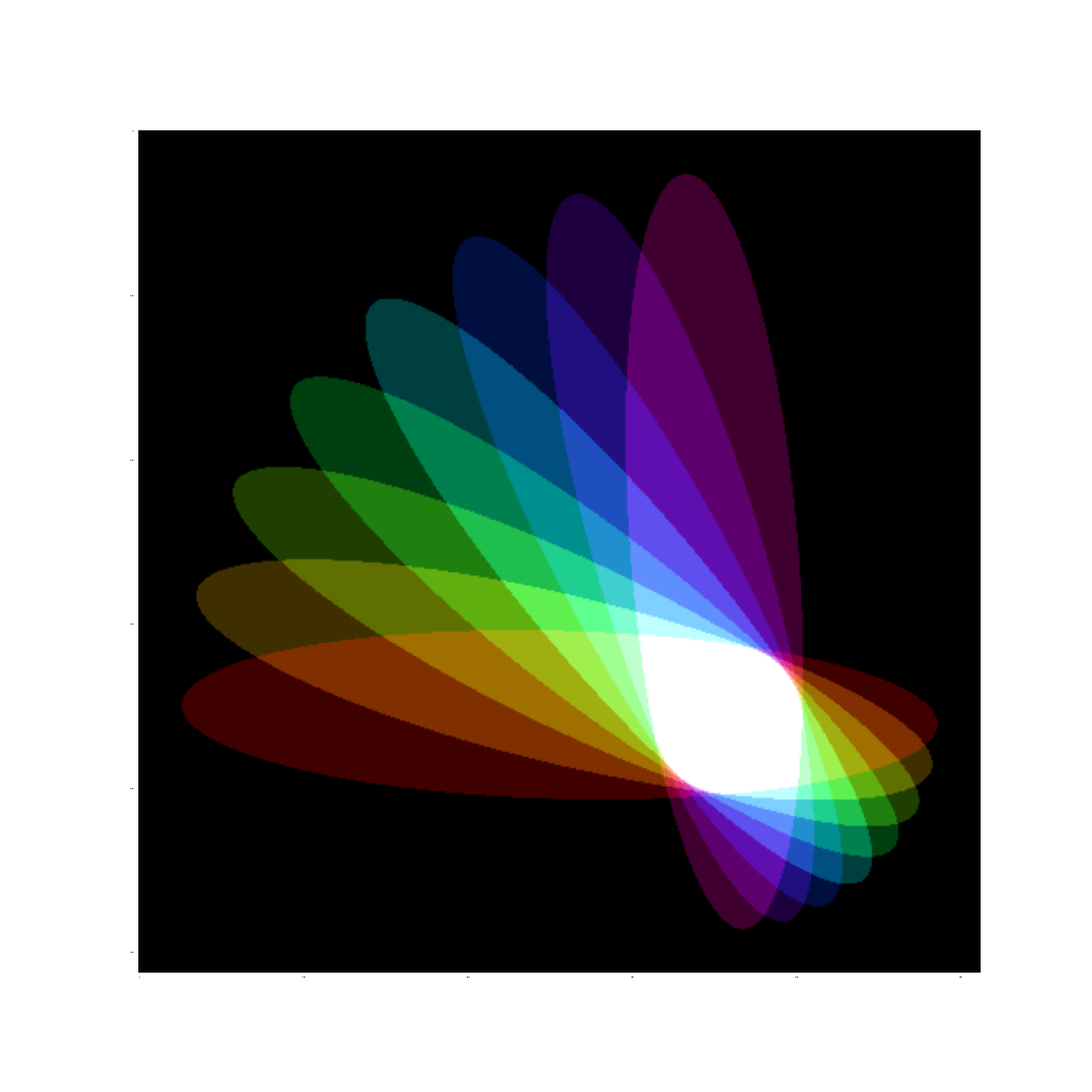}}
	\caption{False colours used for multi-look imagery in \figurename{s}~\ref{fig: multilook}~and~\ref{fig: multilook zoom}. The major axis of each ellipse shown is oriented towards the azimuth of the centre of each sub-aperture.}
\end{wrapfigure}
We apply this to data from the Gotcha carpark challenge dataset\cite{casteel2007challenge}, which is a multiple altitude circular SAR dataset, using the VV-polarised channel of the third pass which is at an elevation angle of approximately \qty{44.7}{\degree}. We take eight subsets of data as a surrogate for multiple ``looks''. These are each \qty{3}{\degree} apertures, with a spacing of \qty{10}{\degree} between each, for a total angular coverage of \qty{84}{\degree}.  A bandwidth of approximately \qty{620}{\mega\hertz} at centre frequency of \qty{9.60}{\giga\hertz} gives each single-channel image a cross-range and range resolution of approximately \SIlist[]{0.2;0.3}{\metre}, respectively.  As this is circular SAR the layover of each subsequent channel will change direction with observation, unlike the case of taking multiple apertures along the same straight flight path (referred to as ``squinted''). Note of course that these effective range and cross-range directions change with each look direction, so the effective resolution in pixel coordinates varies with each image channel.

\begin{figure}
	\centering
	\begin{subfigure}[t]{0.48\textwidth}
		\centering
		\includegraphics[width=\textwidth, trim=1cm 1cm 1cm 1cm, clip]{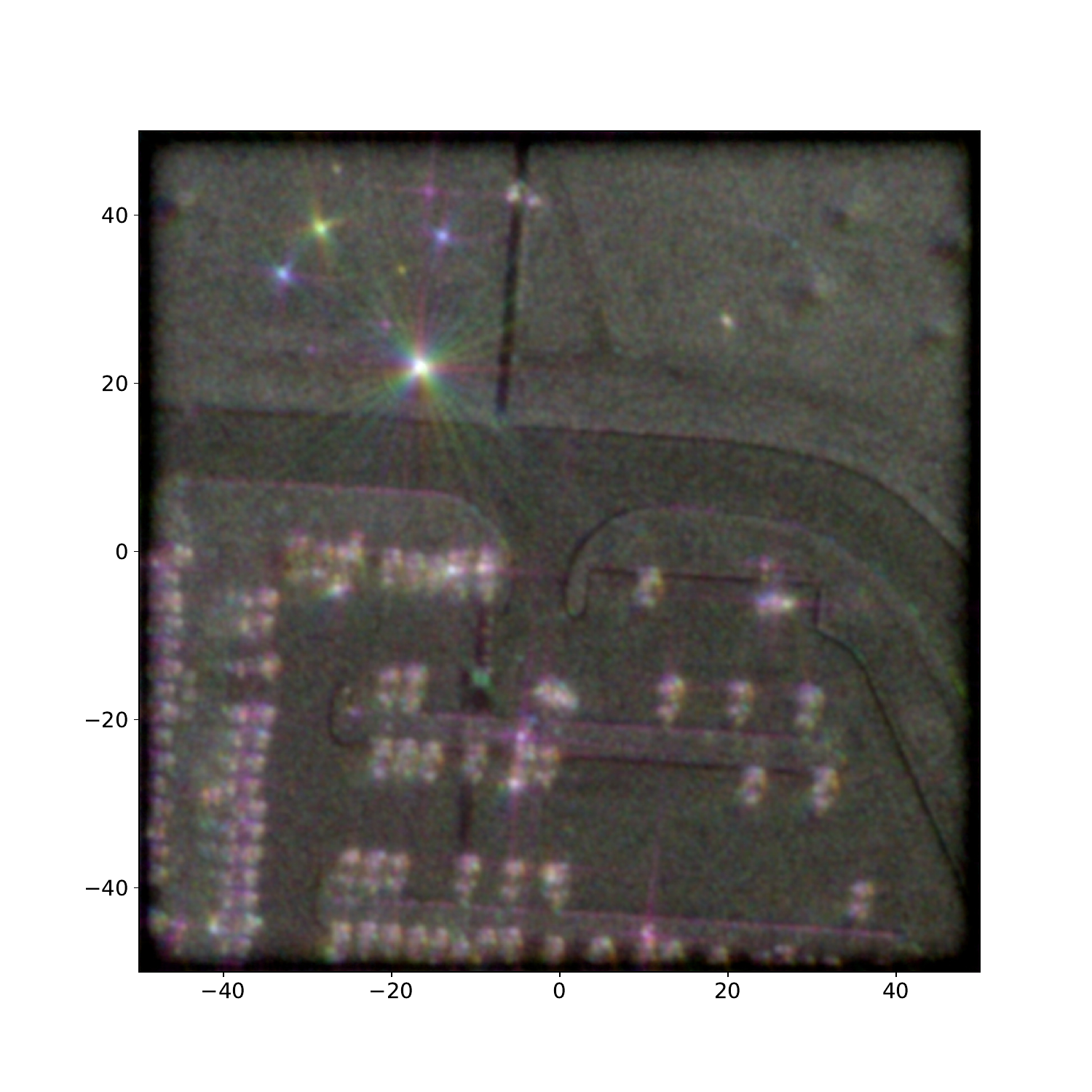}
		\caption{Generalised Tikhonov regularised}
	\end{subfigure}\hfill
	\begin{subfigure}[t]{0.48\textwidth}
		\centering
		\includegraphics[width=\textwidth, trim=1cm 1cm 1cm 1cm, clip]{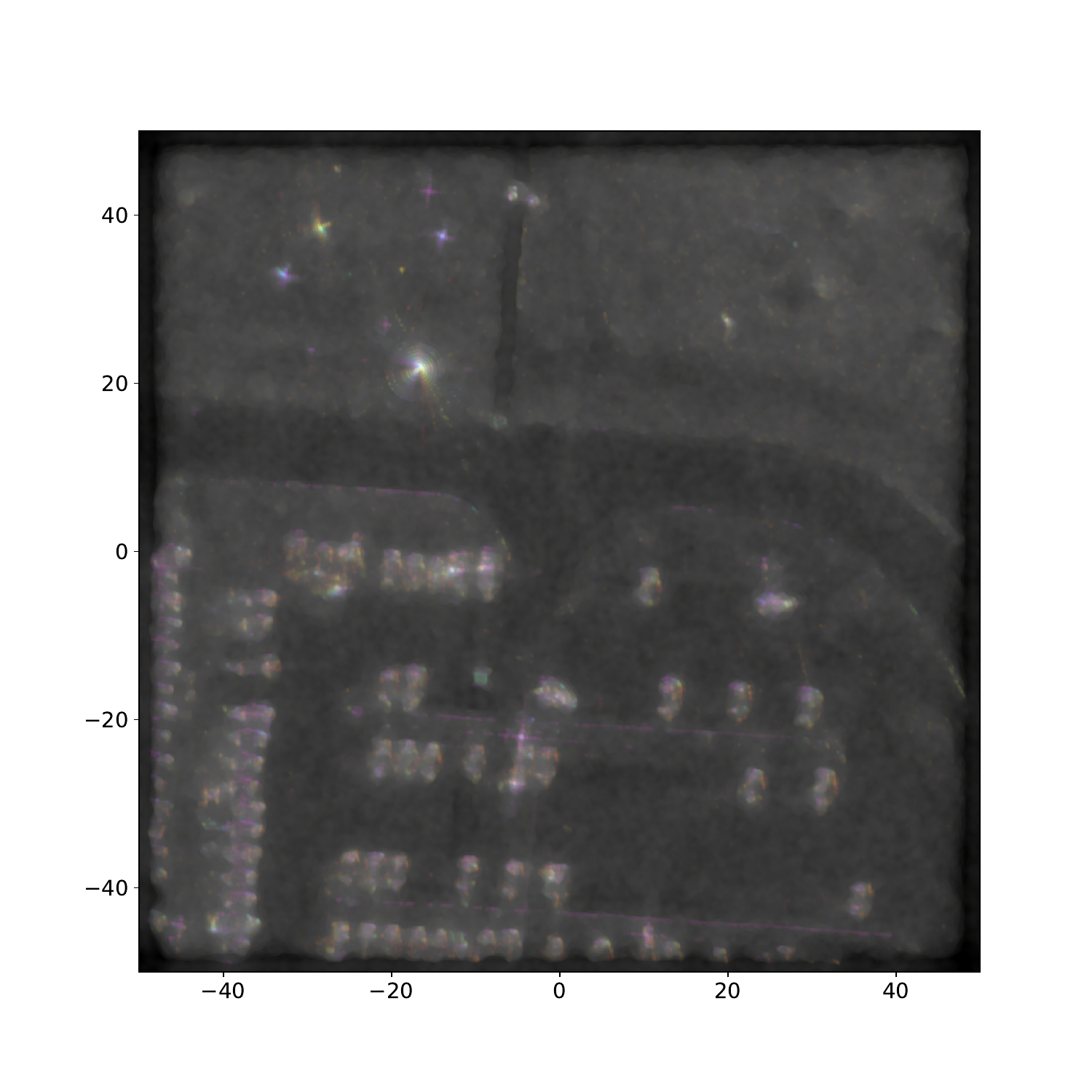}
		\caption{Spatio-channel TV \mbox{regularised}}
	\end{subfigure}~\\
	\begin{subfigure}[t]{0.48\textwidth}
		\centering
		\includegraphics[width=\textwidth, trim=1cm 1cm 1cm 1cm, clip]{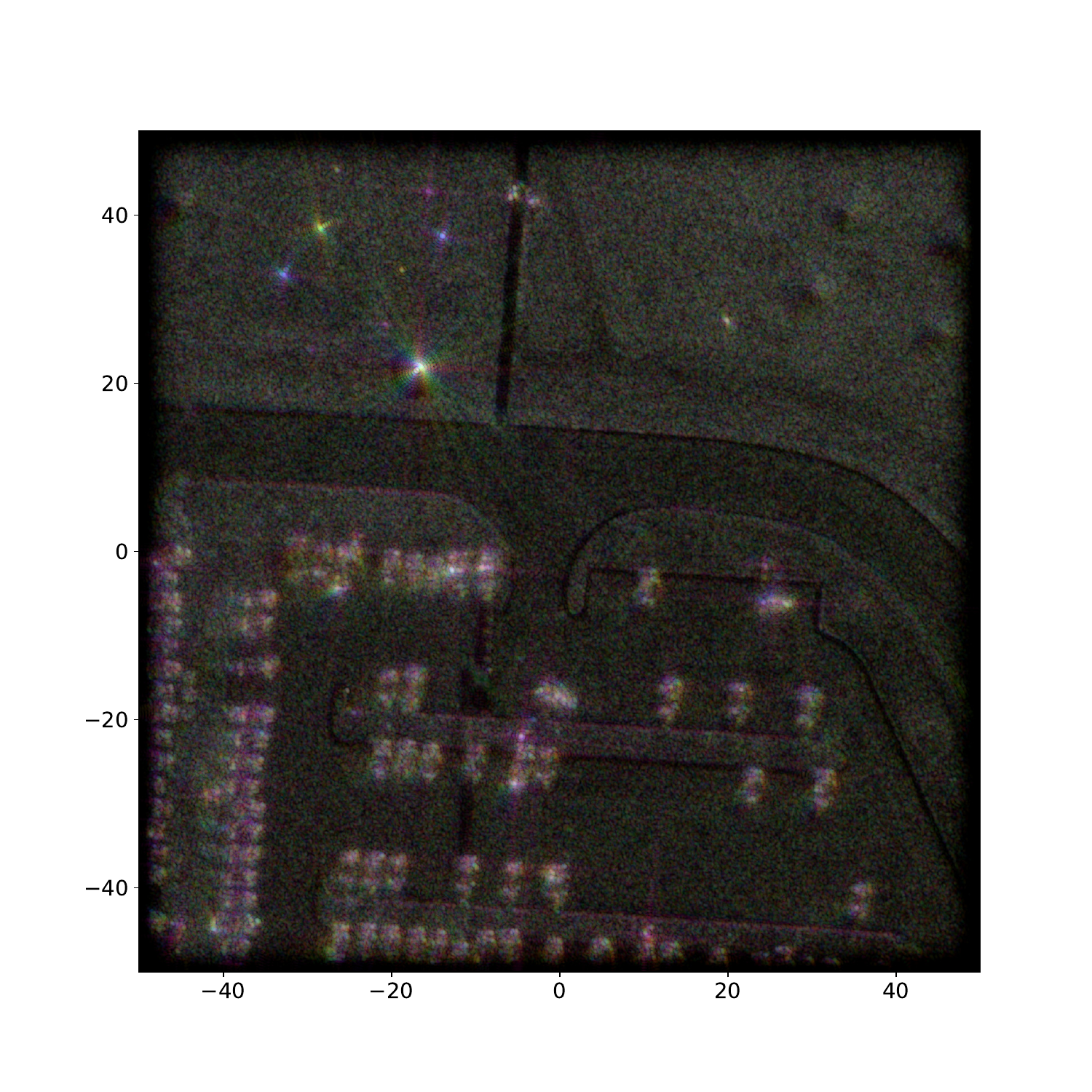}
		\caption{False-colour back-projection}
	\end{subfigure}\hfill
	\begin{subfigure}[t]{0.48\textwidth}
		\centering
		\includegraphics[width=\textwidth, trim=1cm 1cm 1cm 1cm, clip]{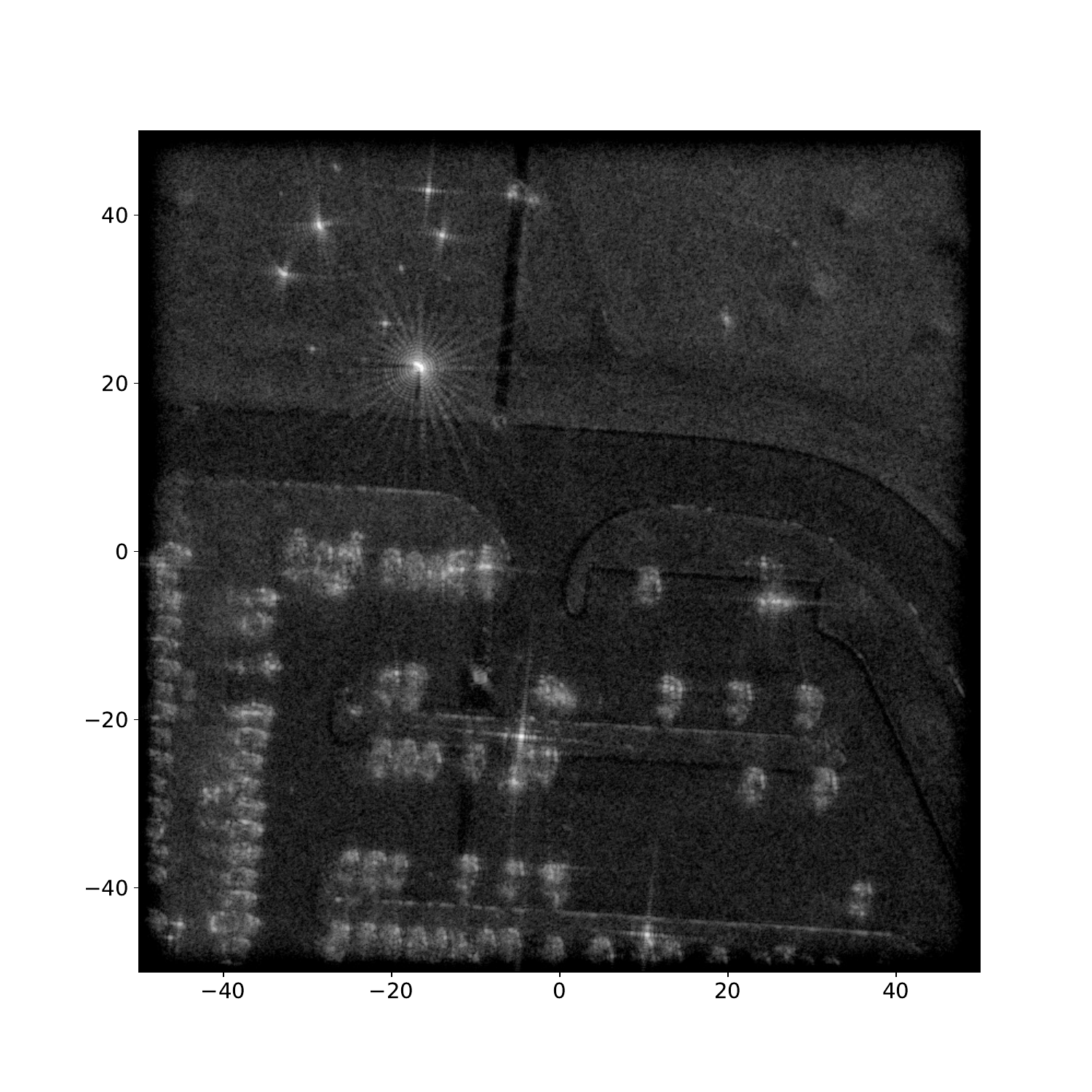}
		\caption{Incoherently combined back-projection}
	\end{subfigure}
	\caption{Multi-aspect images of the Gotcha carpark dataset }
	\label{fig: multilook}
\end{figure}

\begin{figure}
	\centering
	\begin{subfigure}[t]{0.48\textwidth}
		\centering
		\includegraphics[width=\textwidth, trim=1cm 1cm 1cm 1cm, clip]{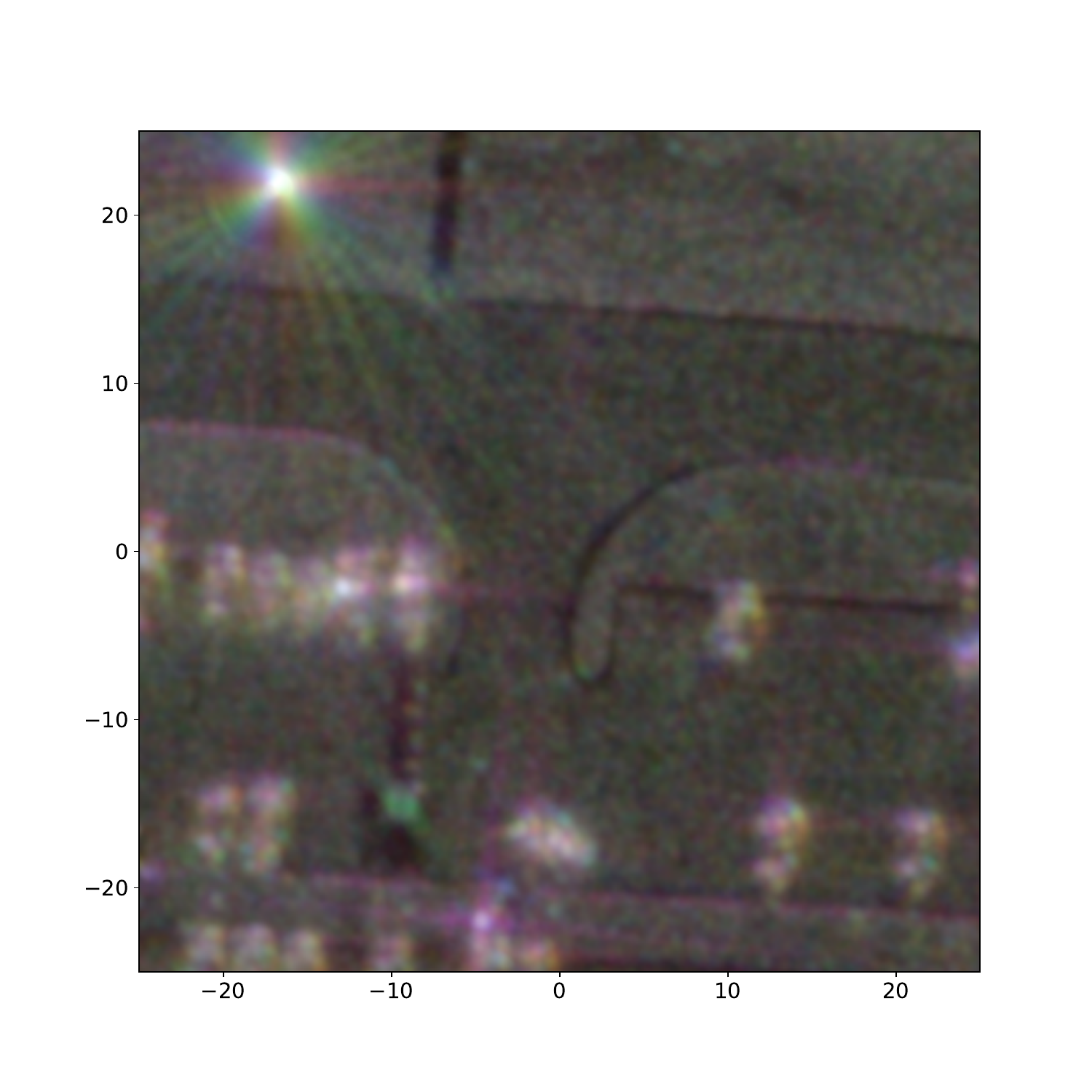}
		\caption{Generalised Tikhonov regularised}
	\end{subfigure}\hfill
	\begin{subfigure}[t]{0.48\textwidth}
		\centering
		\includegraphics[width=\textwidth, trim=1cm 1cm 1cm 1cm, clip]{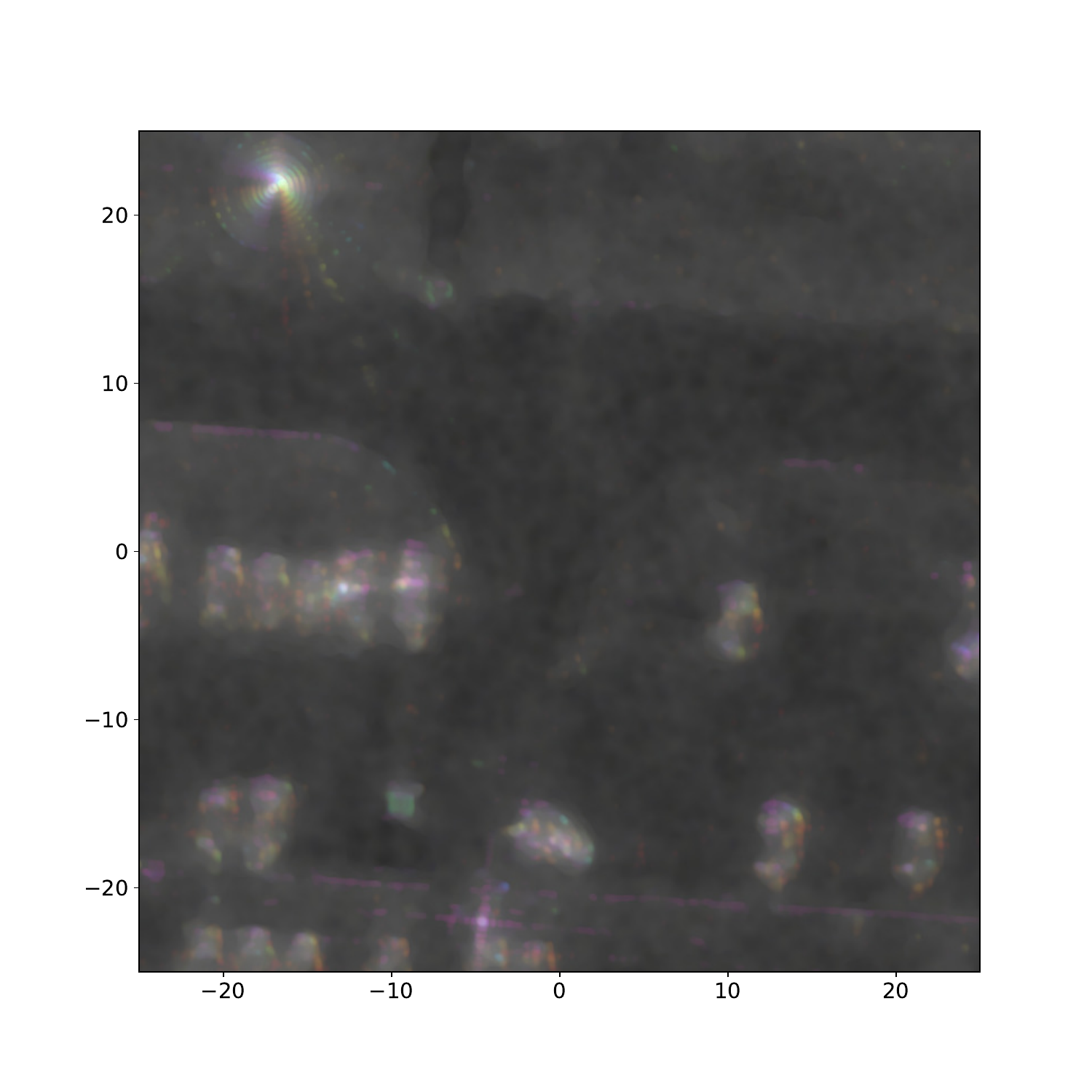}
		\caption{Spatio-channel TV regularised}
	\end{subfigure}~\\
	\begin{subfigure}[t]{0.48\textwidth}
		\centering
		\includegraphics[width=\textwidth, trim=1cm 1cm 1cm 1cm, clip]{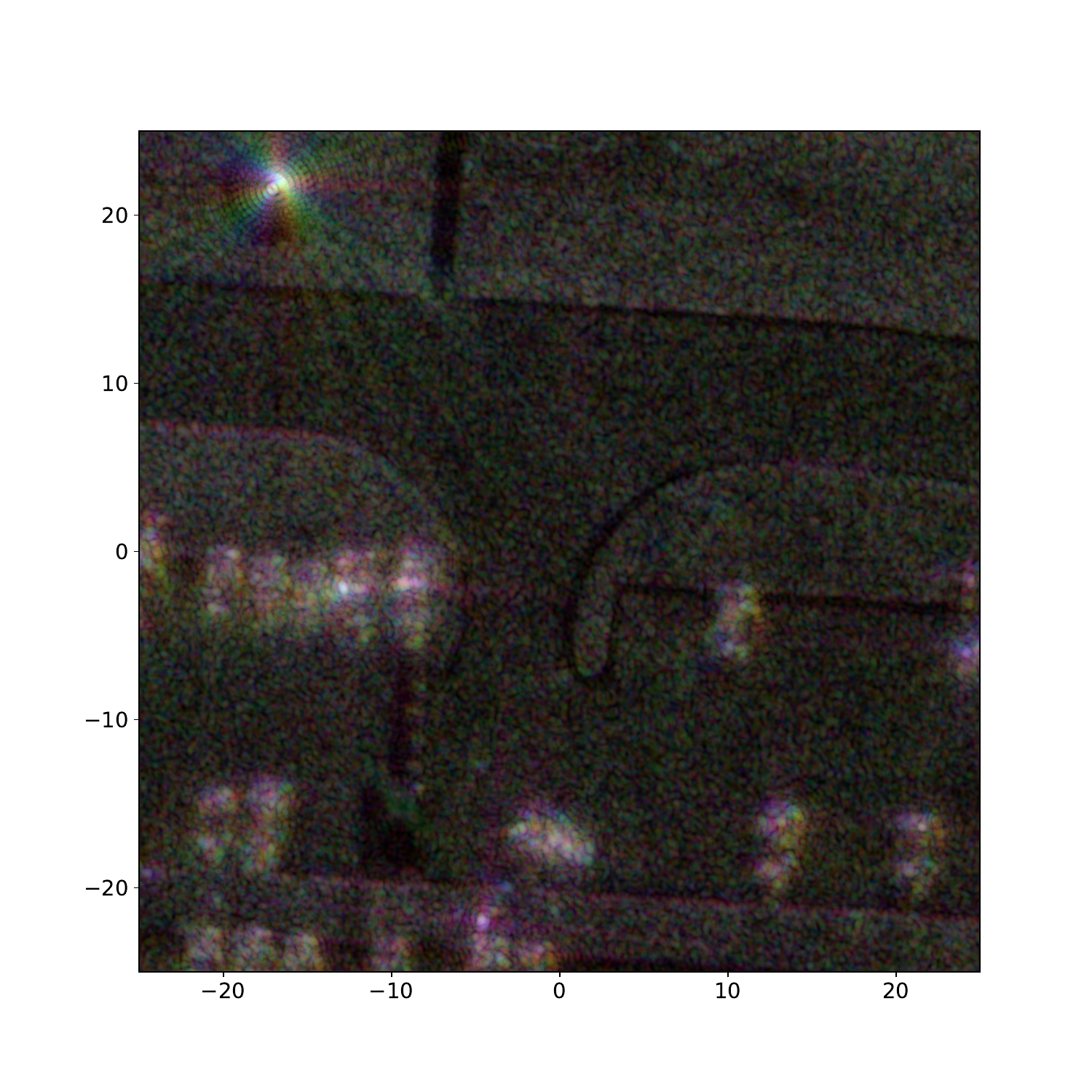}
		\caption{False colour back-projection}
	\end{subfigure}\hfill
	\begin{subfigure}[t]{0.48\textwidth}
	\centering
	\includegraphics[width=\textwidth, trim=1cm 1cm 1cm 1cm, clip]{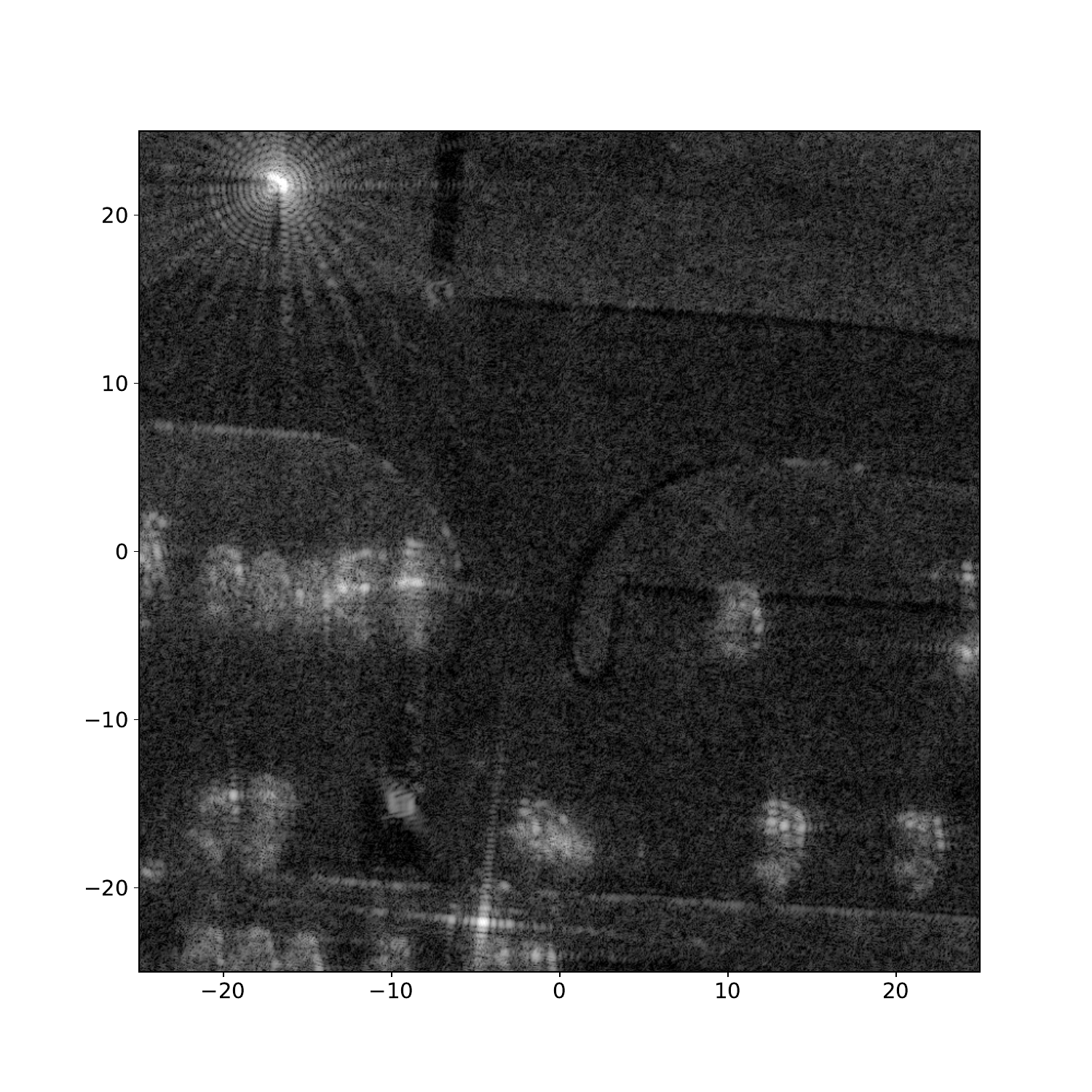}
	\caption{Incoherently combined back-projection}
	\end{subfigure}
	\caption{}
	\label{fig: multilook zoom}
	
\end{figure}

\figurename~\ref{fig: multilook} shows the multi-look reconstruction result, which is also compared to the back-projection as well as the reconstruction result regularised by space-channel TV, described in (\ref{eq: spatiotemporal TV}).  \figurename~\ref{fig: multilook zoom} shows a zoom in of the same results.  These are given a false colour for each channel, such that if a pixel is of equal brightness in each channel then it will appear grey-scale. Each image colour channel is in \unit{\decibel} scale, clipped between \SIrange[range-units=single]{-6}{-31}{\decibel} relative to the peak pixel intensity across all channels.  Each image is formed with a pixel spacing of \qty{0.1}{\metre} in each direction.

We can see a similar speckle reduction in the generalised Tikhonov result as with an incoherent summation image, whilst retaining information about angle-dependent (anisotropic) scattering present similar to the false-colour back-projection.  This combination we find in places also helps to pick out targets appearing more strongly in certain directions, such as the green colourised object at \([-15,-10]\) which is quite weak in the back-projection. It also does something to reduce the dominance of sidelobes from the strong calibration target at \([20,-20]\) seen in the incoherent summation image.  The Space-channel TV regularised image also helps in similar ways, though with a loss of dynamic range (for example obscuring some paths) despite an order of magnitude smaller regularisation parameter.  Purely as a matter of opinion, we find this despeckling and overall appearance of the generalised Tikhonov result pleasing to look at.

While these results are meant only as example reconstructions to highlight one such possible application of Theorem~\ref{thm: main result} to real-world data, we can also make some suggestions as to where a multi-channel generalised Tikhonov regularisation approach may be beneficial.  For example, where multi-aspect imagery has been collected in noisy environments or with strong RF interference, as can be common  in low-frequency SAR systems\cite{doody2017bright}, Tikhonov regularisation may be effective in de-noising and mitigating interference.  Perhaps more so given that any interference is likely to vary in time (and therefore channel/aspect), and so be penalised by differing from the resultant image in neighbouring channels.  There may also be potential for other applications where strong angle-dependent, anisotropic scattering is observed, such as in multi-static data collections.

\FloatBarrier
\subsection{Level set magnitude reconstruction in complex-valued SAR}\label{subsec: level set}
Here we use the same Gotcha carpark dataset to demonstrate the potential of a level set approach applied only to the magnitude of the complex image, as discussed in section~\ref{sec: level set}.  We use the first \qty{3}{\degree} from the second pass, which is at an elevation angle of \qty{45.6}{\degree}. A \qty{499.4}{\mega\hertz} bandwidth subset of the full dataset is taken at the same centre frequency, for approximately equal range and cross-range resolutions of \qty{0.3}{\metre}.

\begin{figure}[h!]
	\centering
	\begin{subfigure}[t]{0.48\textwidth}
		\centering
		\includegraphics[width=0.8\textwidth, trim = 3cm 3.5cm 3cm 3cm, clip]{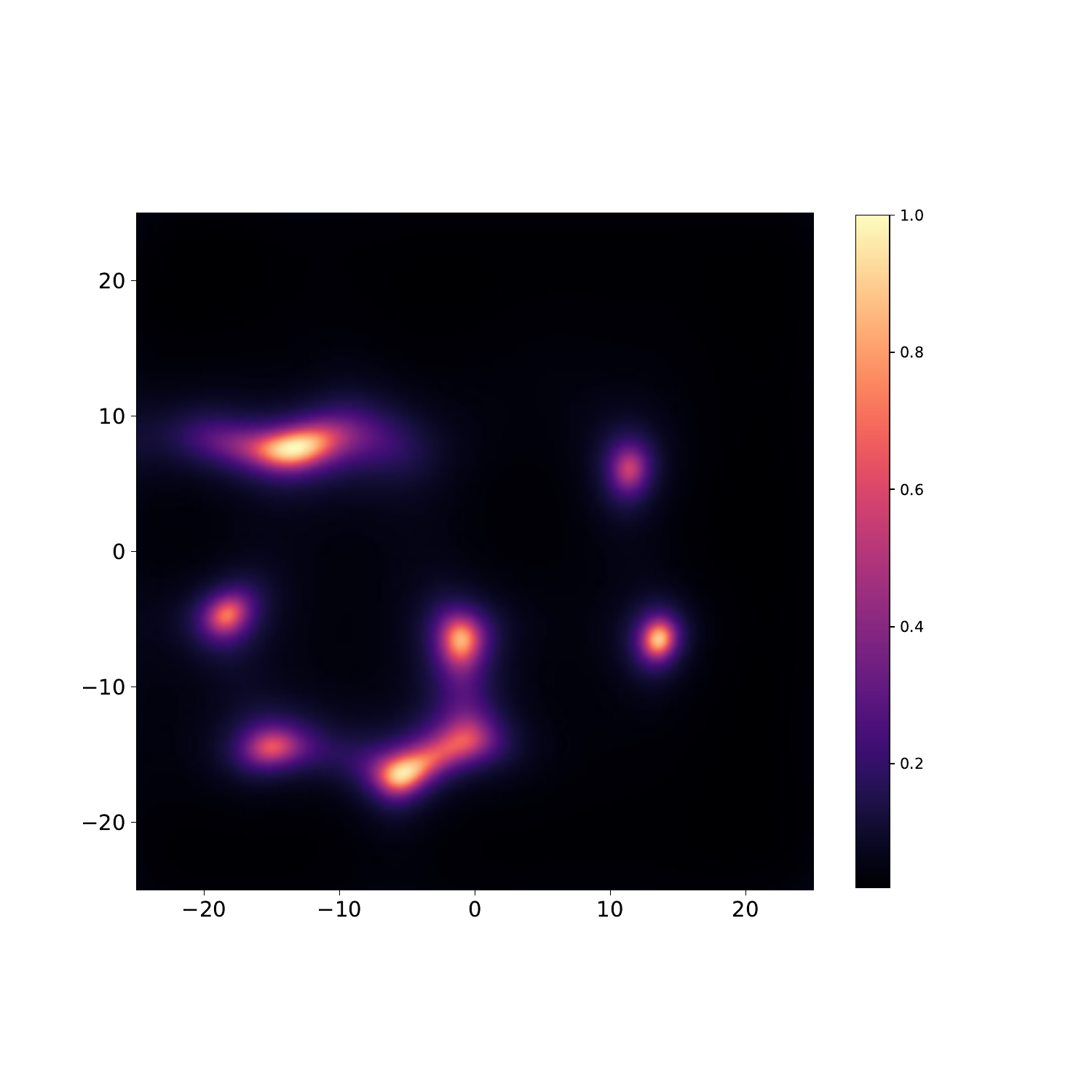}
		\caption{Level set reconstruction, linear scale}
		\label{fig: level lin}
	\end{subfigure}\hfill
	\begin{subfigure}[t]{0.48\textwidth}
		\centering
		\includegraphics[width=0.8\textwidth, trim = 3cm 3.5cm 3cm 3cm, clip]{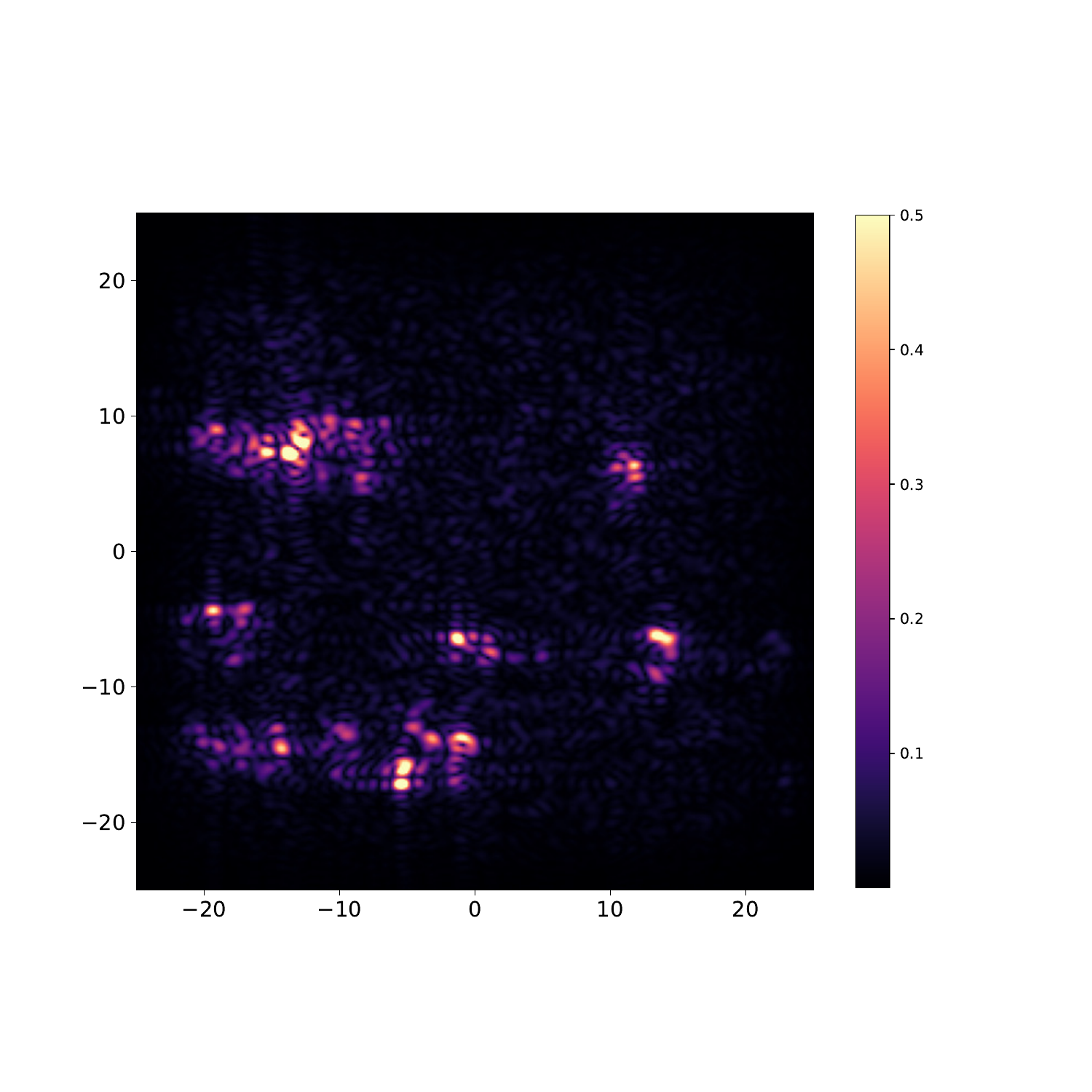}
		\caption{Back-projection, linear scale}
	\end{subfigure}\\
	\begin{subfigure}[t]{0.48\textwidth}
		\centering
		\includegraphics[width=0.8\textwidth, trim = 3cm 3.5cm 3cm 3cm, clip]{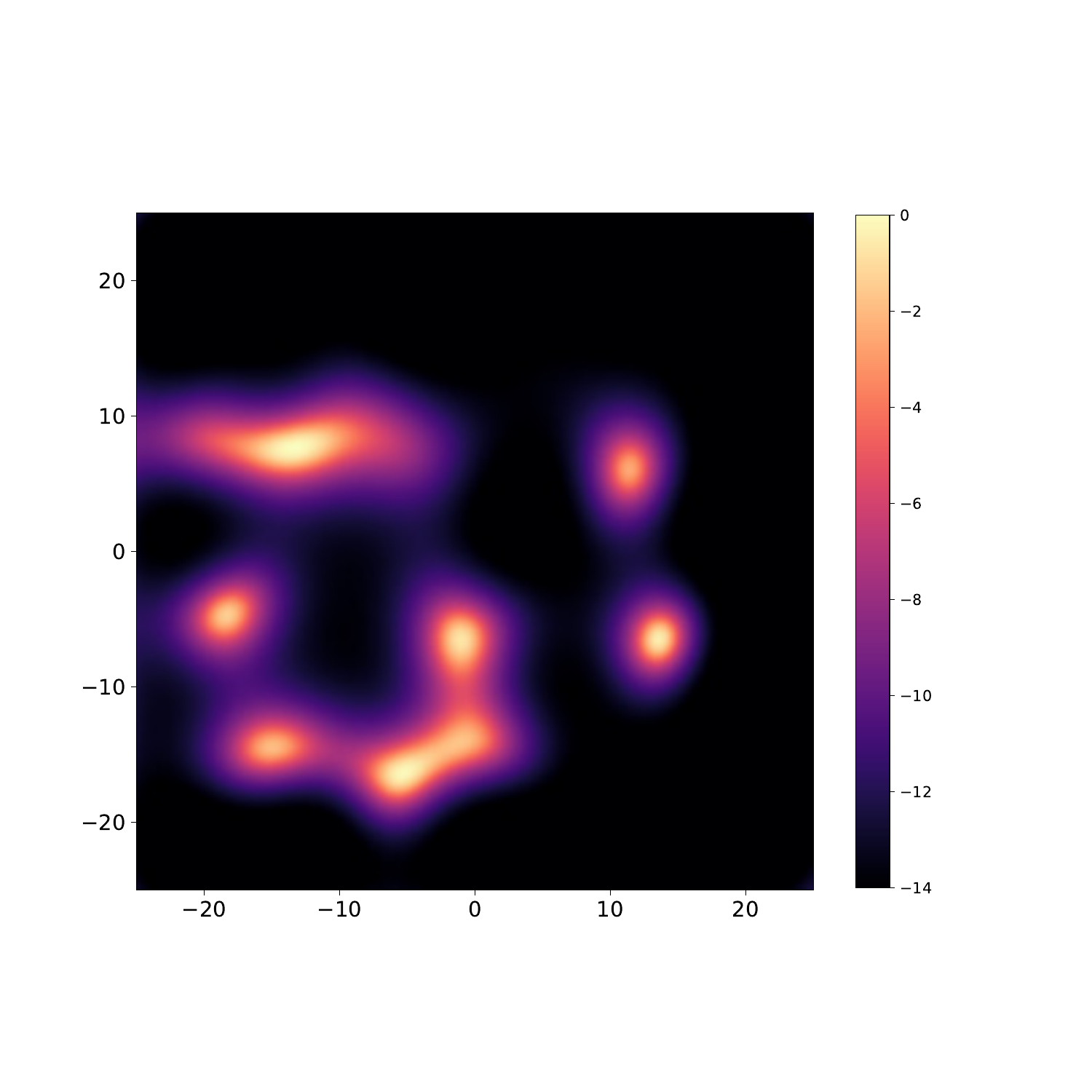}
		\caption{Level set reconstruction, log scale}
	\end{subfigure}\hfill
	\begin{subfigure}[t]{0.48\textwidth}
		\centering
		\includegraphics[width=0.8\textwidth, trim = 3cm 3.5cm 3cm 3cm, clip]{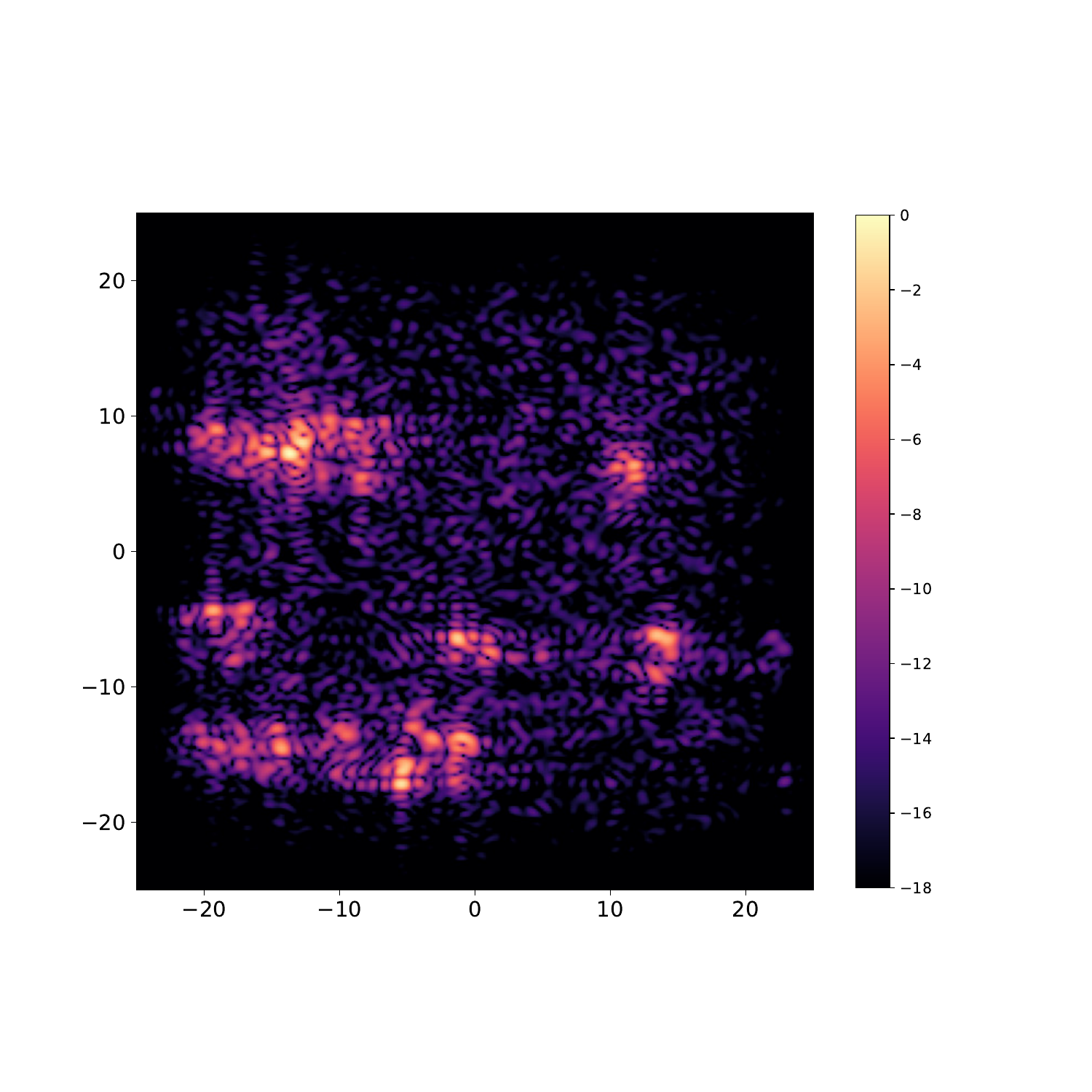}
		\caption{Back-projection, log-scale}
		\label{fig: level bp log}
	\end{subfigure}\\
	\begin{subfigure}[t]{0.48\textwidth}
		\centering
		\includegraphics[width=0.8\textwidth, trim = 3cm 3.5cm 3cm 3cm, clip]{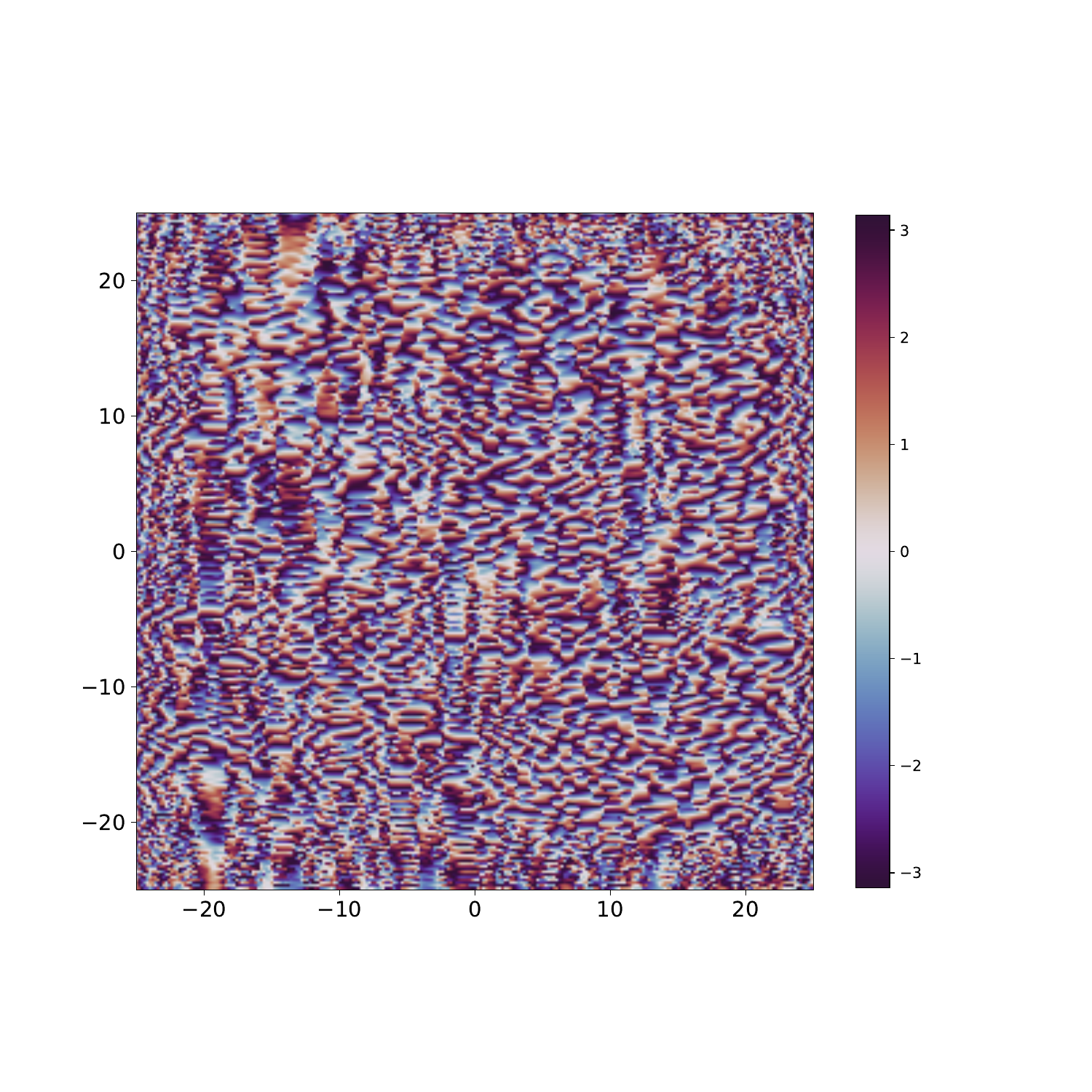}
		\caption{Phase of level set reconstruction}
		\label{fig: level phase}
	\end{subfigure}\hfill
	\begin{subfigure}[t]{0.48\textwidth}
		\centering
		\includegraphics[width=0.8\textwidth, trim = 3cm 3.5cm 3cm 3cm, clip]{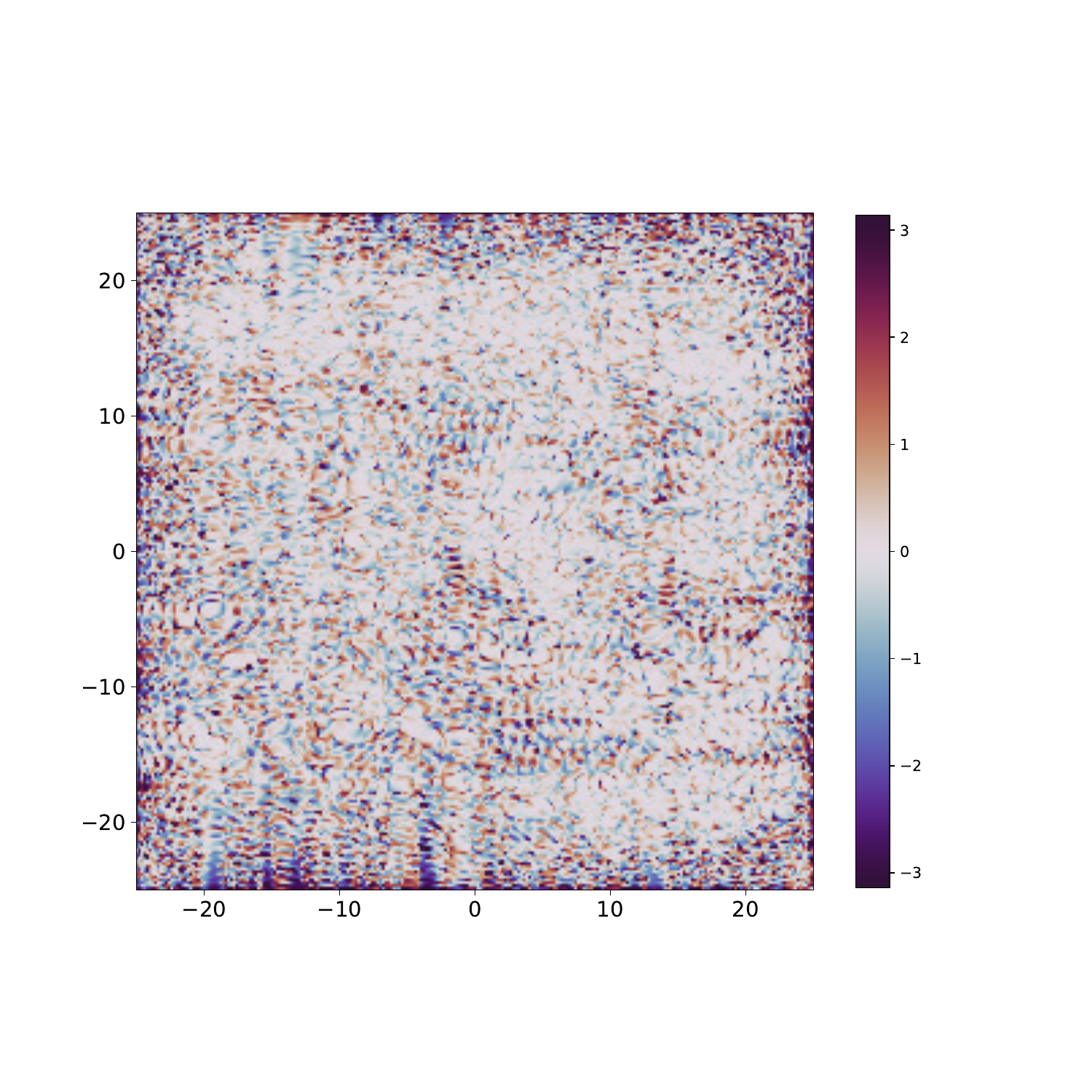}
		\caption{Phase difference between level set reconstruction and backprojection}
		\label{fig: level phase diff}
	\end{subfigure}
	\caption{PaLEnTIR level set reconstruction and back-projection of a sub-scene from the Gotcha dataset. \ref{fig: level lin}-\ref{fig: level bp log} show the magnitude of these complex-valued images, \ref{fig: level phase} the phase of the level set reconstruction, and \ref{fig: level phase diff} the phase difference between level set reconstruction and backprojection images}
	\label{fig: palentir gotcha}
\end{figure}

Reconstructing only a central \qty{50}{\meter\squared} region, the result using PaLEnTIR\cite{ozsar2024parametric} is shown in \figurename~\ref{fig: palentir gotcha} and compared to the back-projection.  We can see the parts of vehicles are brought out against the background speckle. Since this particular level-set formulation is expressive in terms of smoothly varying edges of objects, the specular nature of SAR is also homogenised in a smoothly-varying manner in the level set reconstruction. From \figurename~\ref{fig: level phase} we can see this reconstruction does indeed vary pixel-to-pixel in phase, and from \figurename~\ref{fig: level phase diff} we see the resulting phases are generally close to those of the backprojection though not identically so.

Such reconstructions could be useful in pre-screening, detection and classification of targets against relatively strong clutter, noise and interference.  One might even be able to use the coefficients of the radial basis functions themselves as a means to classify detected objects.

As with the discussion in section~\ref{sec: level set}, we propose the method here as a basis for more computationally efficient methods to apply level sets only to the magnitude of complex-valued imagery. For example, we have found some computational benefit in under-solving the proximal map (\ref{eq: prox level set}) in earlier iterations.  This suggests one potential approach generalising (\ref{eq: level set image denoise}) to change the indicator function of the level set basis \(\chi_{\mathcal{C}}\) to a smoothly-varying metric of the distance from \(\Span\{f(\vec{p})\}\) (for example based on smooth approximations of the Heaviside step function).  It may also be possible to include Jacobian information or weighting into the projection step to speed convergence.  Alternative level set methods themselves (i.e. other than PaLEnTIR) may also be useful for SAR in this framework.

\subsection{Total Generalised Variation for spotlight mode SAR}\label{subsec: tgv}
Total Generalised Variation (TGV) extends the notion of Total Variation to one which balances the first \(k\) derivatives of a function, which reduces the so-called staircasing effect\cite{bredies2010total}.  The second order total generalised variation function is written in discrete form as
\begin{equation}
	\TGV_{\alpha\beta}^2(\lvec{u}) = \min_{\lvec{w}} \alpha\|\lmat{D}\lvec{u}-\lvec{w}\|_{2,1} + \beta\|\mathcal{E}\lvec{w}\|_{2,1},
	\label{eq: TGV}
\end{equation}
where \(\mathcal{E}\) is the discrete symmetrised gradient operator and \(\lmat{D}\) the discrete gradient operator. We can see that using \(\TGV_{\alpha,\beta}^2\) as a regularisation term will promote images with piecewise-constant-gradient. Moreover, as \(\alpha/\beta\rightarrow 0\), (\ref{eq: TGV}) tends to the (isotropic) Total Variation.  \(\TGV_{\alpha,\beta}^2\) has proximal map
\begin{equation}
	\begin{aligned}
		\tilde{\lvec{u}} =& \prox_{\TGV_{\alpha,\beta}^2}(\lvec{y}),\\
		(\tilde{\lvec{u}}, \tilde{\lvec{w}}) =& \argmin_{\lvec{u},\lvec{w}}\frac{1}{2}\|\lvec{u}-\lvec{y}\|_2^2 + \alpha\|\lmat{D}\lvec{u}-\lvec{w}\|_{2,1} + \beta\|\mathcal{E}\lvec{w}\|_{2,1},
	\end{aligned}
	\label{eq: TGV prox}
\end{equation}
which may be solved numerically for example via TGV\cite{papoutsellis2021core2}.

Papafitsoros and Bredies have shown exact and numerical solutions to (\ref{eq: TGV prox}) for which \(\tilde{\lvec{u}}=\prox_{\TGV_{\alpha,\beta}^2}(\lvec{r})\notin\mathds{R}^n_{\geq0}\) with \(\lvec{r}\in\mathds{R}^n_{\geq0}\) in the case of 1D functions\cite{papafitsoros2015study}. From the structure of these solutions, we can assume that the solution to \(\prox_{\TGV_{\alpha,\beta}^2}^+(\lvec{r})\) will neither be as simple as taking the element-wise absolute value of \(\tilde{\lvec{u}}\), nor projecting into \(\mathds{R}^n_{\geq0}\), since both would increase the second derivative at places where \(\tilde{\lvec{u}}\) crosses an axis.  While a positivity constraint may be incorporated into the numerical solution to (\ref{eq: TGV prox}) directly, this may not always be possible (for example if a black-box code has been used), but Algorithm~\ref{alg: prox abs} may always be applied.

We demonstrate the application of Algorithm~\ref{alg: prox abs} to TGV-regularised least-squares reconstructions, i.e.
\begin{equation}
	\minimise \|\lmat{A}\lvec{z}-\lvec{d}\|_2^2 + \lambda \TGV_{\alpha,\beta}^2(|\lvec{z}|),
	\label{eq: TGV regularised}
\end{equation}
using CPHD format data from the Umbra Open Data Program\cite{umbra}. The sub-scene is of the Diamond Light Source. A centre frequency of \qty{9.8}{\giga\hertz}, bandwidth of \qty{672}{\mega\hertz}, and synthetic aperture length of \qty{2.5}{\degree}, provides a range and cross-range resolution of approximately \SIlist[list-units=single]{0.22;0.35}{\metre} in the slant plane, respectively (giving approximately square pixels in the ground plane).

For an example reconstruction, we take a regularisation parameter of \(\lambda=10\) which we find ensures the effect of TGV is clear in the resulting images.  Within each of these outer optimisation iterations, calculating \(\prox_{\TGV_{\alpha,\beta}^2}\), i.e. solving (\ref{eq: TGV prox}), is also carried out via PDHG and is allowed 100 iterations for this (much less computationally expensive) inner optimisation problem. A thorough test of appropriate stopping conditions is beyond the scope of this paper, but we generally observe the solution has stagnated before this point.

\begin{figure}[h]
	\centering
	\begin{subfigure}[t]{\textwidth}
		\centering
		\includegraphics[width=0.8\textwidth, trim = 0.9cm 3cm 2cm 5cm, clip]{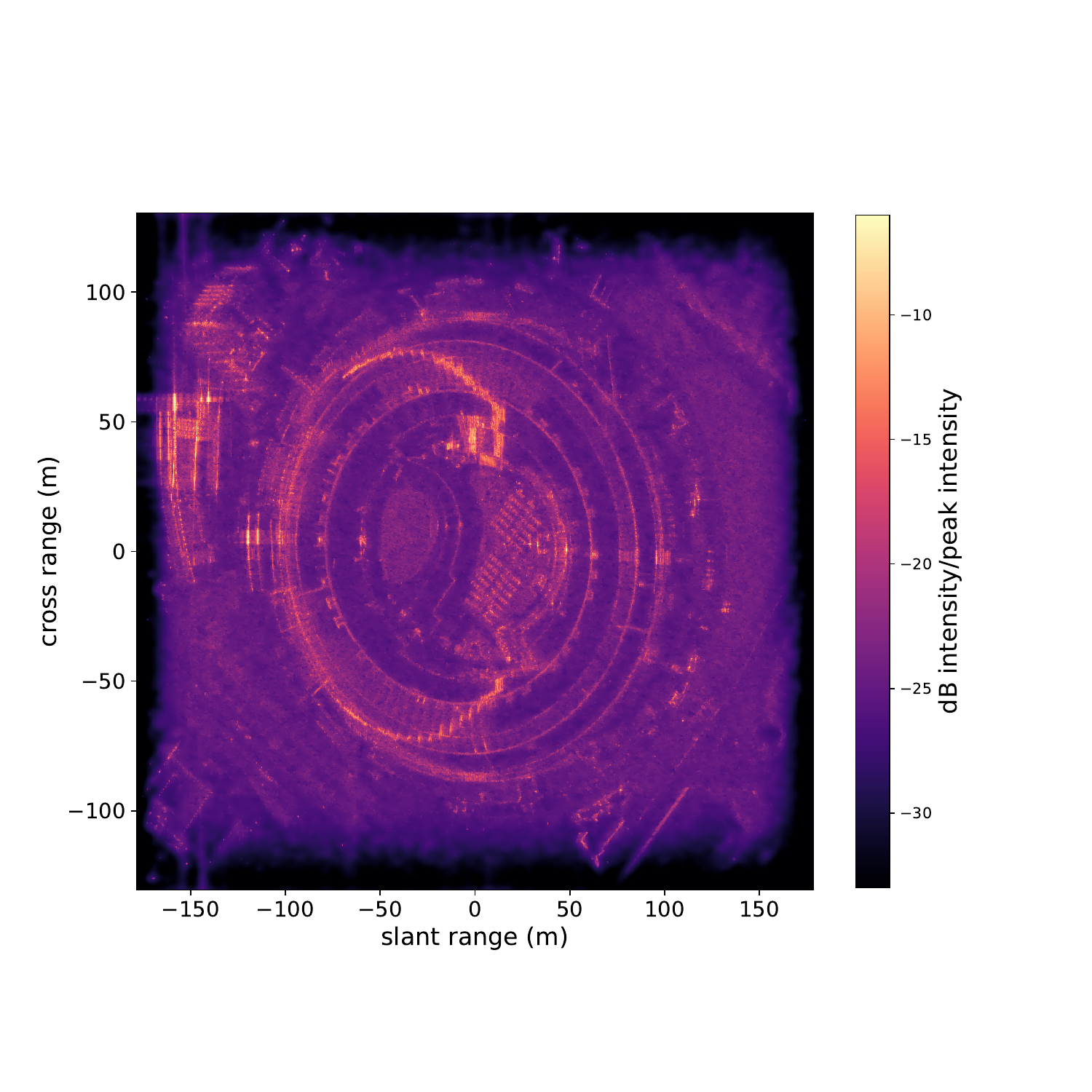}
		\caption{\(\TGV_{1,2}^2\)-regularised reconstruction}
	\end{subfigure}\\
	\begin{subfigure}[t]{\textwidth}
		\centering
		\includegraphics[width=0.8\textwidth, trim = 0.9cm 3cm 2cm 3.5cm, clip]{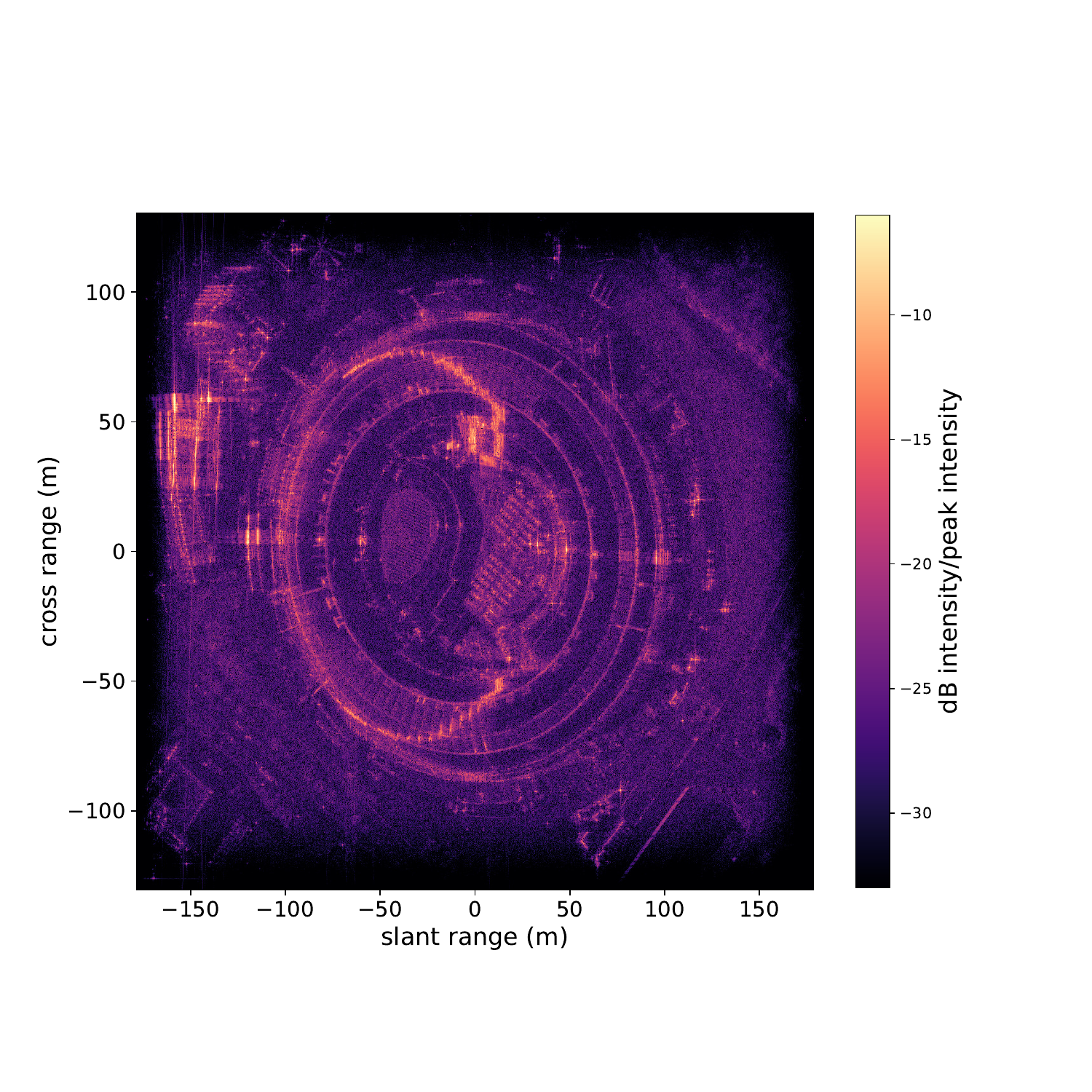}
		\caption{Back-projection}
		\label{fig: diamond full}
	\end{subfigure}\hfill
	
\caption{TGV-regularised reconstruction and back-projection images of Umbra CPHD SAR data\cite{umbra} of the Diamond Light Source.  These have \(2041\times1677\) pixels.}
\end{figure}

\begin{figure}[h]
\centering
	\begin{subfigure}[t]{0.48\textwidth}
		\centering
		\includegraphics[width=\textwidth, trim = 1.cm 3cm 2cm 5cm, clip]{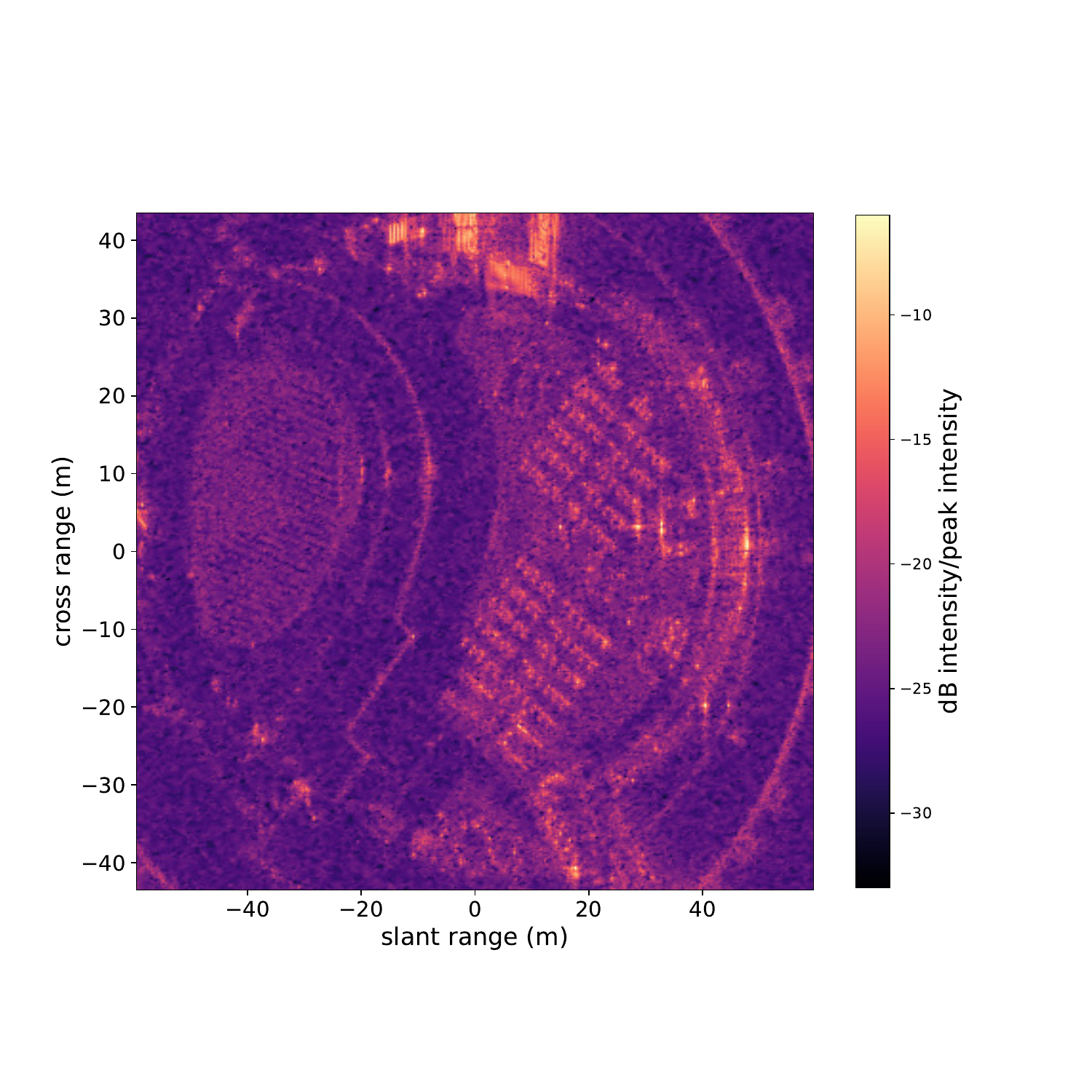}
		\caption{\(\TGV_{1,0.25}^2\)}
	\end{subfigure}\hfill
	\begin{subfigure}[t]{0.48\textwidth}
		\centering
		\includegraphics[width=\textwidth, trim = 1cm 3cm 2cm 4cm, clip]{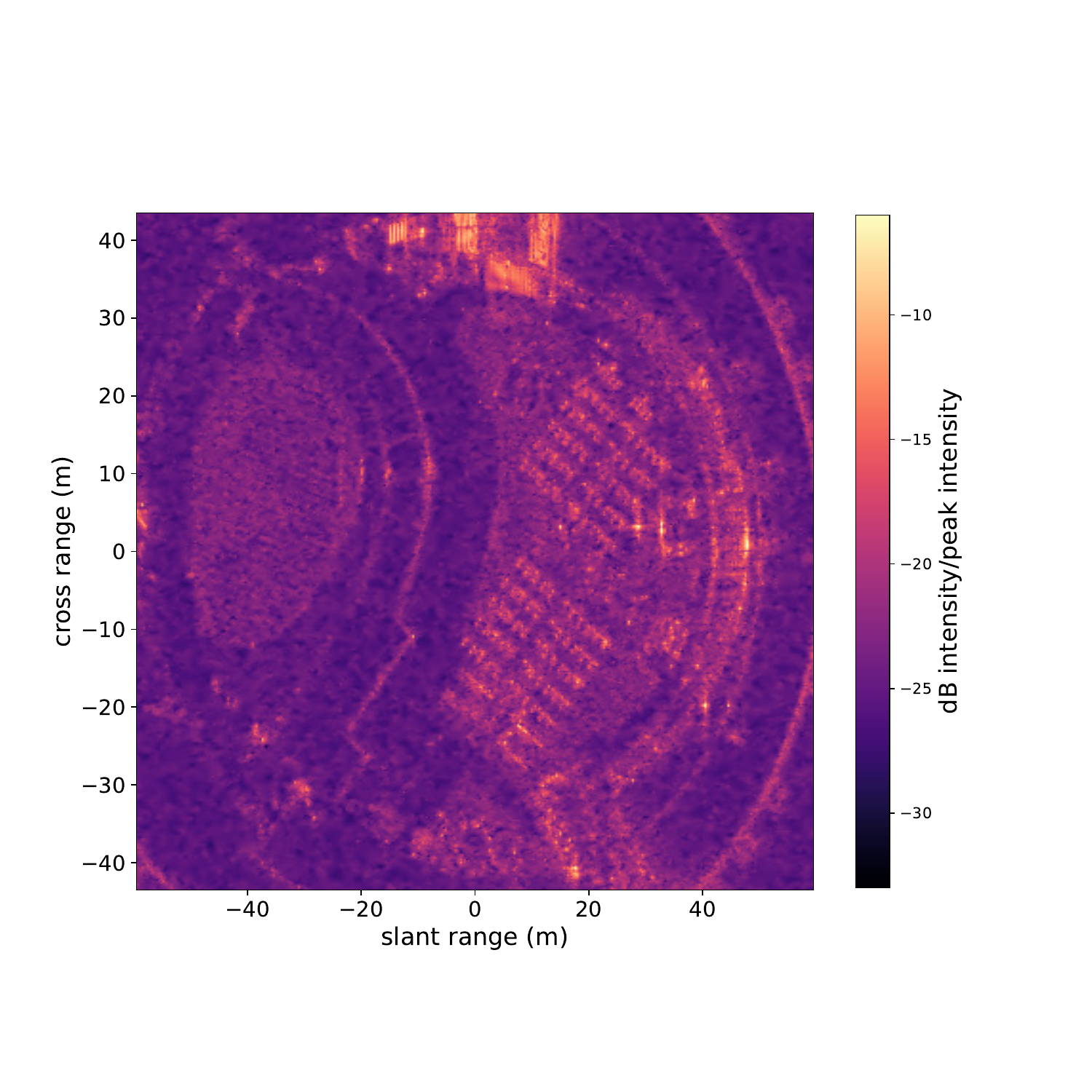}
		\caption{\(\TGV_{1,0.5}^2\)}
	\end{subfigure}\\
	\begin{subfigure}[t]{0.48\textwidth}
		\centering
		\includegraphics[width=\textwidth, trim = 1cm 3cm 2cm 5cm, clip]{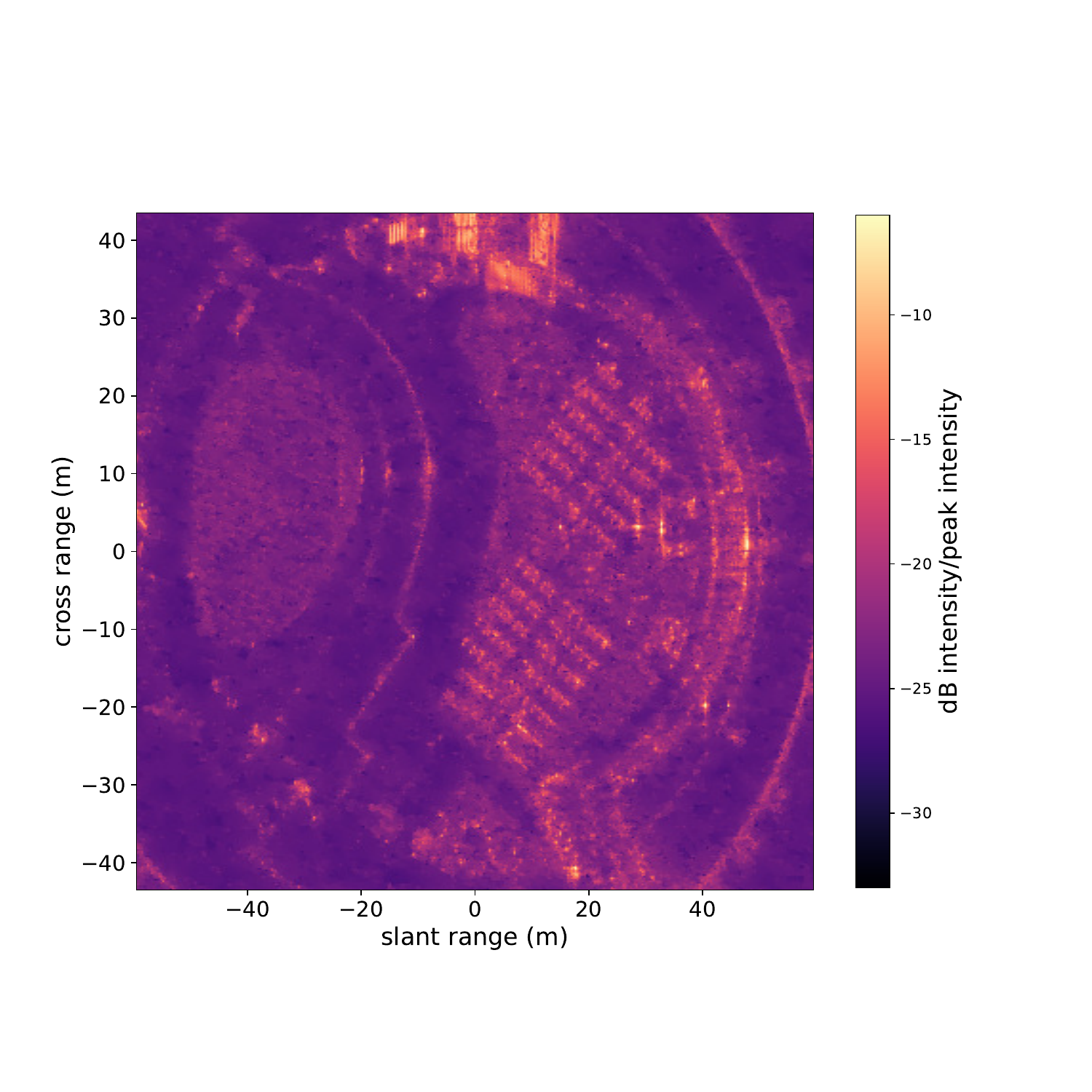}
		\caption{\(\TGV_{1,2}^2\)}
	\end{subfigure}\hfill
	\begin{subfigure}[t]{0.48\textwidth}
		\centering
		\includegraphics[width=\textwidth, trim = 1cm 3cm 2cm 4cm, clip]{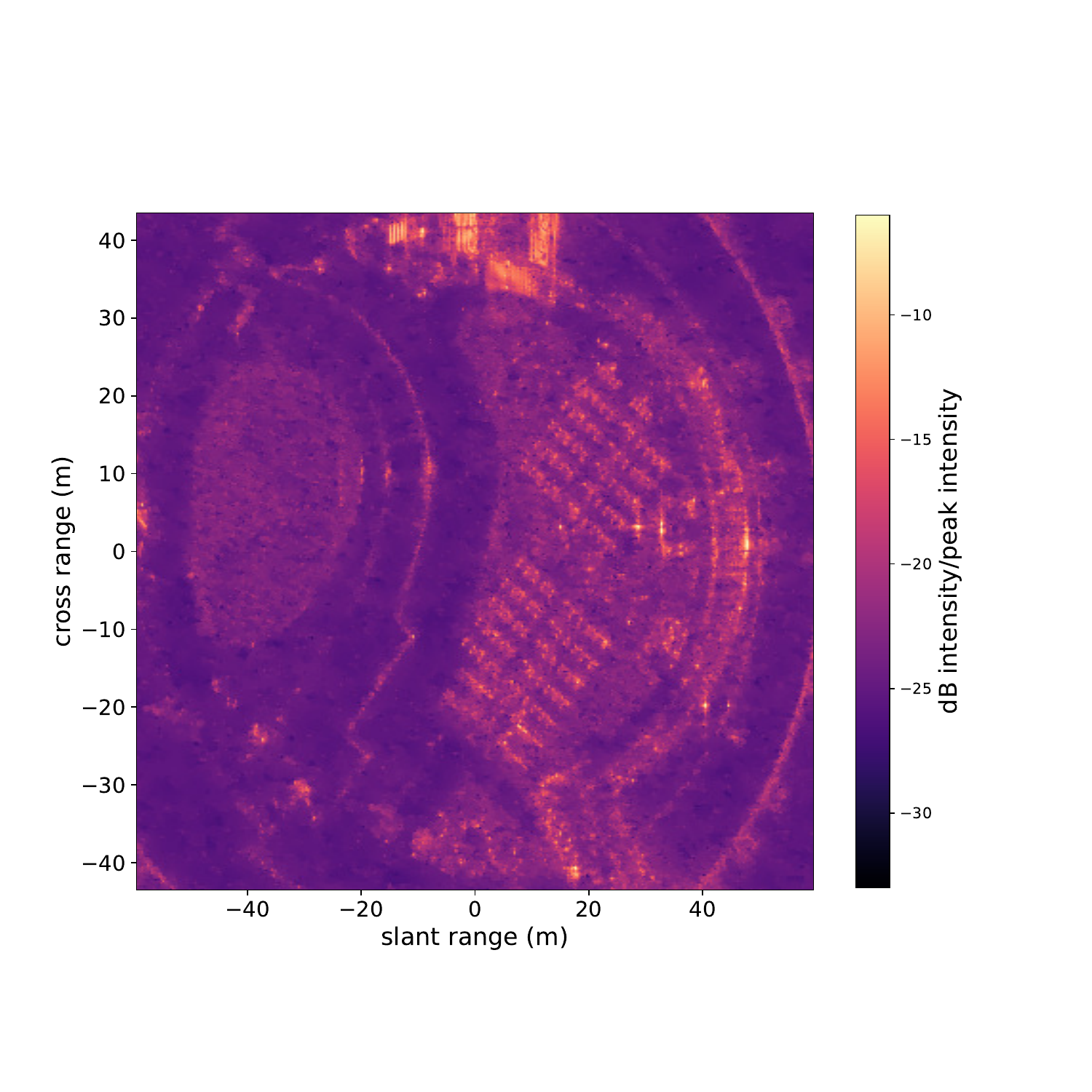}
		\caption{\(\TGV_{1,8}^2\)}
		\label{subfig: TGV18 diamond}
	\end{subfigure}\\
	\begin{subfigure}[t]{0.48\textwidth}
	\centering
	\includegraphics[width=\textwidth, trim = 1cm 3cm 2cm 5cm, clip]{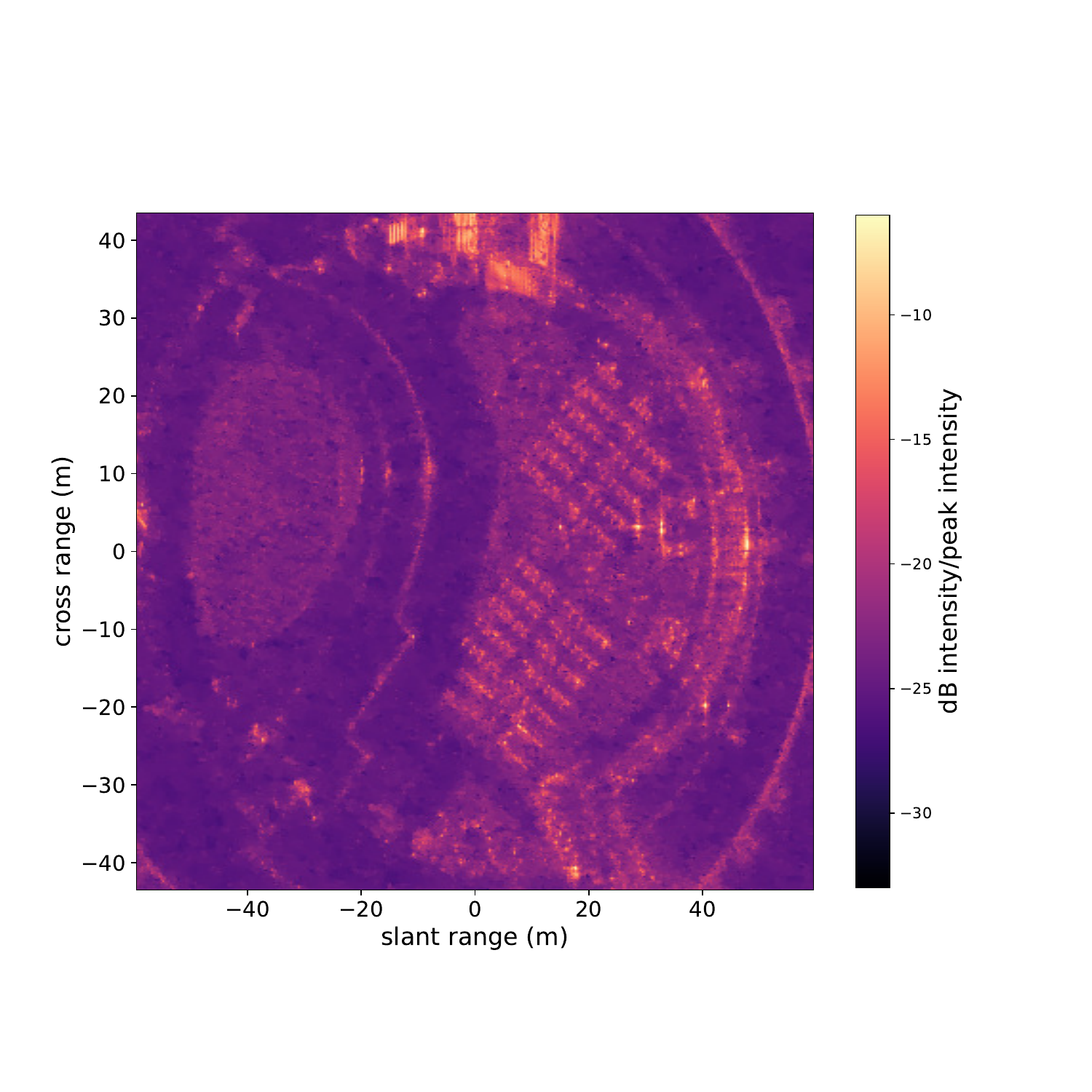}
	\caption{\(\TV\equiv\TGV_{1,\infty}^2\)}
	\label{subfig: TV diamond}
	\end{subfigure}\hfill
	\begin{subfigure}[t]{0.48\textwidth}
	\centering
	\includegraphics[width=\textwidth, trim = 1cm 3cm 2cm 4cm, clip]{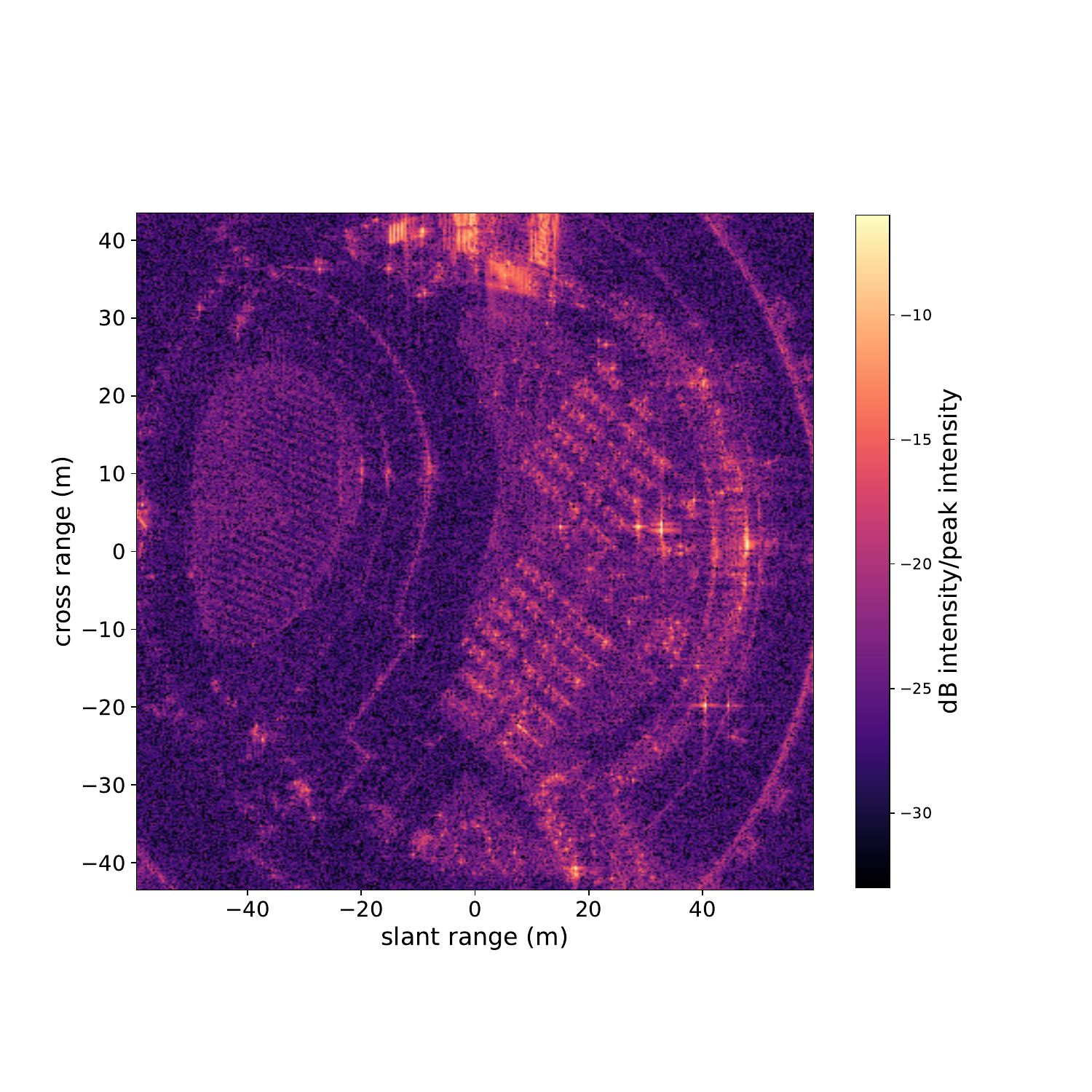}
	\caption{Back-projection}
	\end{subfigure}	\caption{Comparison of TGV and TV-regularised reconstruction, and back-projection images of Diamond Light Source CPHD SAR data via Umbra open data program\cite{umbra}, showing central zoomed in image chips.}
	\label{fig: diamond zoom}
\end{figure}

\figurename~\ref{fig: diamond full} shows the \(\TGV_{1,2}^2\)-regularised reconstruction and back-projection of the full extent of data used.  
\figurename~\ref{fig: diamond zoom} shows zoomed image chips to compare the TGV-regularised reconstructions for \(\TGV_{1,0.25}^2\), \(\TGV_{1,0.5}^2\), \(\TGV_{1,2}^2\) and \(\TGV_{1,8}^2\), as well as the \(\TV\)-regularised and back-projection images.  The TGV reconstructions have significant speckle reduction versus the back-projection image, appear to have some sidelobe suppression of strong scatterers, whilst also appearing to maintain the structure of objects in the scene. For smaller \(\beta\) values, finer low-contrast structural features are retained, but speckle begins to be reintroduced. As mentioned above, taking \(\beta/\alpha\rightarrow\infty\) results in an image which is more like a TV reconstruction, i.e. piecewise constant. Comparing \figurename~\ref{subfig: TGV18 diamond} with \figurename~\ref{subfig: TV diamond}, we might conclude that \(\infty\approx8\) in so far as the visual effect of TGV regularisation mirrors that of TV in this particular case.

For each of the \(\TGV_{\alpha\beta}^2\), we found that only the first few iterations ever needed to enter the Douglas-Rachford iterations of Algorithm~\ref{alg: prox abs} to calculate a bounded proximal map.  This is perhaps due to the highly speckled nature of standard SAR (backprojection) imagery, with high valued pixels adjacent to those of negligible strength.  The earlier iterations will share these features (starting with a step in the backprojection direction), resulting in the initial \(\prox_{\TGV}\) step of Algorithm~\ref{alg: prox abs} crossing below zero similar to the aforementioned results of Papafitsoros \textit{et al}\cite{papafitsoros2015study}.  After a few PDHG iterations, the TGV measure of the iterate has decreased, speckle has been smoothed out, and so the \(\prox_{\TGV}\) step no longer attempts to cross out of the positive orthant to compensate for high frequency, high contrast pixel-to-pixel changes.  Thus, Algorithm~\ref{alg: prox abs} only adds some small additional computational cost in this instance.

\section{Conclusion}
We have considered the problem of applying regularisation \(H\) to the magnitude only of complex-valued (coherent) reconstruction problems, \(H(|\cdot|)\), with a particular motivation being synthetic aperture radar.  Many optimisation algorithms used to solve general reconstruction algorithms make use of the proximal map, \(\prox_H\), particularly where \(H\) is non-smooth.  We have shown how one can simply calculate this proximal map of \(H(|\cdot|)\) applied only to the magnitude of the complex image under a certain (broad) sufficient condition, i.e. \(\prox_{H(|\cdot|)}\), namely that the proximal map will always lie in the positive orthant, without needing to solve a non-linear phase fitting problem. This makes use of the proximal map of \(H\) itself.  We demonstrate that several commonly used regularisation functions meet this sufficient condition, and so can be readily applied to the magnitude of complex imagery.  Moreover, the main result is used to provide an initial, simple means to apply level set reconstruction methods to the magnitude of coherent reconstruction problems, though we later discuss that further development may be needed to make this computationally practical.

Furthermore, we provide a simple algorithm to calculate the proximal map of \(H(|\cdot|)\) for other functions which may not meet the sufficient conditions of our main result.  Since this does not introduce any additional computational cost where the main result would hold, it may be applied in a black-box manner by users without consideration as to whether specific theoretical conditions hold.  This is beneficial to users of coherent imagery, for whom we have made fairly arbitrary regularisation functions readily available.  To demonstrate the real-world applicability of our results, we provide example reconstructions using publicly available real airborne and satellite SAR data.  For these we employ generalised Tikhonov, level sets, and total generalised variation, to the magnitude of both multi-look and single-channel SAR imagery, which to our knowledge has not previously appeared in the literature.

\FloatBarrier

\ack
The authors would like to thank the Isaac Newton Institute for Mathematical Sciences, Cambridge, for support and hospitality during the programme ``Rich and Nonlinear Tomography: A Multidisciplinary Approach'' (RNT), where work on this paper was undertaken. This work was supported by EPSRC grant no EP/R014604/1 and EP/V007742/1, and has made use of computational support by CoSeC, the Computational Science Centre for Research Communities, through CCPi.  Watson was supported by the Royal Academy of Engineering and the Office of the Chief Science Adviser for National Security under the UK Intelligence Community Postdoctoral Research Fellowship programme. 

We would also like to thank Misha Kilmer for useful discussions on level set reconstruction methods, Eric Miller and Ege Ozsar for provision of PaLEnTIR code, as well as Evangelos Papoutsellis for useful discussions on total generalised variation.

\FloatBarrier

\begin{appendices}
\section{SAR data model}
\label{ap: data model}
The modelling of SAR data are based on a single-scattering assumption (the Born approximation) of scalar waves from a stationary reflectivity function \(V\) \cite{watson2024resolving, cheney2009fundamentals}. That is, in the frequency domain,
\begin{equation}
	U^{\mathrm{sc}}(\svec{x},\omega) =-\int G_0(\svec{x}-\svec{z},\omega)V(z)\omega^2U^{\mathrm{in}}(\omega,\svec{z}) \rmd \svec{z},
\end{equation}
where \(G_0\) is the free-space Green's function for the Helmholtz equation, 
\begin{equation}
	G_0(\svec{x},\omega)=\frac{\rme^{\imi \omega|\svec{x}|}}{4\pi|\svec{x}|},
\end{equation}
with \(c\) the free-space wave speed, and \(U^{\mathrm{in}}\) is the incident wavefield.  For a point-like antenna, this can be modelled as a multiple of the freespace Green's function itself,
\begin{equation}
	U^{\mathrm{in}}(\svec{z},\omega) = p(\omega)G_0(\svec{z}-\svec{\gamma}_T),
\end{equation}
where \(\svec{\gamma}_T\) is the position of radar transmitter.  If the scene is modelled by isotropic point scatterers, \(V(\svec{z}):=\sum_i v_i\delta(\svec{z}_i-\svec{z})\), then (writing out the Green's functions explicitly) the scattered field at receiver location \(\svec{\gamma}_R\) is given as
\begin{equation}
	U^{\mathrm{sc}}(\svec{\gamma}_R,\omega) = a(\omega)\sum_i v_i \rme^{i\omega(|\svec{\gamma_R}-\svec{z}_i|+|\svec{z}_i-\svec{\gamma}_T|)/c},
\end{equation}
where the \(\omega^2\) and amplitude terms have been absorbed into \(a(\omega):=\omega^2P(\omega)/(4\pi r)\). Moreover, we have assumed the loss due to distance from the scene to the radar is approximately constant. This is commonly assumed at in SAR processing, though can easily be included where necessary such as imaging at shorter ranges.  

The measured data itself will be the scattered field I/Q demodulated by mixing with a reference signal, which will be the reflection expected from a reference location \(\svec{x}_{\mathrm{ref}}\) (the \textit{scene reference point}).  If the transmit and receiving antennas move along (possibly co-located) paths \(\svec{\gamma}_T(s)\), \(\svec{\gamma}_R(s)\) parameterised by \(s\) (referred to as \textit{slow-time}), then we can write
\begin{equation}
	d(s,\omega) = U^{\mathrm{sc}}(s,\omega)\rme^{\imi\omega(|\svec{\gamma_R}-\svec{x}_{\mathrm{ref}}|+|\svec{x}_{\mathrm{ref}}-\svec{\gamma}_T|)/c}.
	\label{eq: data continuous d}
\end{equation}
This assumes platform motion is negligible during the time-of-flight of the pulse (the \textit{start-stop} approximation).  Since pulses will be transmitted at a discrete set of times, and digitised and recorded at a discrete set of frequencies, we can write the data model discretely as
\begin{equation}
	\lvec{d}=\lmat{A}\lvec{v}, \quad \lvec{v}=[v_1,\ldots,v_N], \quad \lvec{d}=[d(s_1,\omega_1), \ldots ,d(s_n, \omega_m)].
	\label{eq: data discrete d}
\end{equation}

The assumptions of single scattering of scalar waves from stationary, isotropic point scatterers are almost universally made in using SAR data, and are often reasonable.  Importantly, they allow for fast and efficient SAR image formation algorithms based on applying the adjoint of (\ref{eq: data continuous d}) (or an approximation to it) -- i.e. back-projection.

It is interesting to note that it is the combination of these assumptions as well as the band-limited data which results in a complex-valued image in which the phase has little meaning for a single image taken in isolation.  In particular, the assumption of scattering due to isotropic point scatterers at pixel locations chosen by us: such locations may not be the exact location (and range) of a true scatterer -- if a single scatterer in within the resolution cell is indeed a good description.  The result is that phase of the image is highly sensitive to small positional errors in a way in which the magnitude of the response is not.

The data model (\ref{eq: data continuous d} - \ref{eq: data discrete d}) may be
used directly in a least-squares reconstruction, but faster implementations
are possible in the time-domain by making use of efficient FFT algorithms.
Time delays between emission, scattering, and reception of impulses are calculated,
and then projected to or from an upsampled, $0^\mathrm{th}$ order-interpolated
delay profile, in forward or adjoint mode respectively.
As the upsampling becomes finer, the time-domain calculation
result approaches the frequency domain result and so we can achieve
arbitrary accuracy.
Forward and adjoint evaluations in the frequency domain have a time
complexity scaling as $\mathcal{O}(Nnm)$,
as we must sum over every combination of scattering point, slow time and
fast frequency. In the time domain, the fast frequency dependence essentially
drops out of the time complexity
at the cost of an insignificant increase memory usage, leaving us with
$\mathcal{O}(Nn)$, a significant improvement as problems grow larger.
When using the time-domain model, we still store data (and calculate the data misfit)
in the frequency domain for memory management purposes, since the cost of the FFTs
is small.

The adjoint form of this time-domain implementation is commonly referred to as the ``fast back-projection algorithm'' in the SAR literature\cite{doerry2016bp}.

\end{appendices}

\printbibliography

@article{bredies2010total,
	title={Total generalized variation},
	author={Bredies, Kristian and Kunisch, Karl and Pock, Thomas},
	journal={SIAM Journal on Imaging Sciences},
	volume={3},
	number={3},
	pages={492--526},
	year={2010},
	publisher={SIAM}
}

@article{aghamiry2021complex,
	title={Complex-valued imaging with total variation regularization: an application to full-waveform inversion in visco-acoustic media},
	author={Aghamiry, Hossein S and Gholami, Ali and Operto, Stephane},
	journal={SIAM Journal on Imaging Sciences},
	volume={14},
	number={1},
	pages={58--91},
	year={2021},
	publisher={SIAM}
}

@ARTICLE{zhao2012separate,
	author={Zhao, Feng and Noll, Douglas C. and Nielsen, Jon-Fredrik and Fessler, Jeffrey A.},
	journal={IEEE Transactions on Medical Imaging}, 
	title={Separate Magnitude and Phase Regularization via Compressed Sensing}, 
	year={2012},
	volume={31},
	number={9},
	pages={1713-1723},
	doi={10.1109/TMI.2012.2196707}
}

@ARTICLE{cetin2001feature,
	author={Cetin, M. and Karl, W.C.},
	journal={IEEE Transactions on Image Processing}, 
	title={Feature-enhanced synthetic aperture radar image formation based on nonquadratic regularization}, 
	year={2001},
	volume={10},
	number={4},
	pages={623-631},	
	doi={10.1109/83.913596}
}

@article{guven2016augmented,
	title={An augmented Lagrangian method for complex-valued compressed SAR imaging},
	author={G{\"u}ven, H Emre and G{\"u}ng{\"o}r, Alper and Cetin, M{\"u}jdat},
	journal={IEEE Transactions on Computational Imaging},
	volume={2},
	number={3},
	pages={235--250},
	year={2016},
	publisher={IEEE}
}

@article{chambolle2004algorithm,
	title={An algorithm for total variation minimization and applications},
	author={Chambolle, Antonin},
	journal={Journal of Mathematical imaging and vision},
	volume={20},
	pages={89--97},
	year={2004},
	publisher={Springer}
}

@article{holman2020emission,
	title={Emission tomography with a multi-bang assumption on attenuation},
	author={Holman, Sean and Richardson, Philip},
	journal={arXiv preprint arXiv:2001.04190},
	year={2020}
}

@article{rambour2019introducing,
	title={Introducing spatial regularization in SAR tomography reconstruction},
	author={Rambour, Cl{\'e}ment and Denis, Lo{\"\i}c and Tupin, Florence and Oriot, H{\'e}l{\`e}ne M},
	journal={IEEE Transactions on Geoscience and Remote Sensing},
	volume={57},
	number={11},
	pages={8600--8617},
	year={2019},
	publisher={IEEE}
}

@article{marques2011sar,
	title={SAR image segmentation based on level set approach and $\{$$\backslash$cal G$\}$ \_A\^{} 0 model},
	author={Marques, Regis C Pinheiro and Medeiros, F{\'a}tima N and Nobre, Juvencio Santos},
	journal={IEEE transactions on pattern analysis and machine intelligence},
	volume={34},
	number={10},
	pages={2046--2057},
	year={2011},
	publisher={IEEE}
}

@article{osher2001level,
	title={Level set methods: an overview and some recent results},
	author={Osher, Stanley and Fedkiw, Ronald P},
	journal={Journal of Computational physics},
	volume={169},
	number={2},
	pages={463--502},
	year={2001},
	publisher={Elsevier}
}

@misc{ozsar2024parametric,
	title={Parametric Level-sets Enhanced To Improve Reconstruction (PaLEnTIR)}, 
	author={Ege Ozsar and Misha Kilmer and Eric Miller and Eric de Sturler and Arvind Saibaba},
	year={2024},
	eprint={2204.09815},
	archivePrefix={arXiv}
}

@book{cheney2009fundamentals,
	title={Fundamentals of radar imaging},
	author={Cheney, Margaret and Borden, Brett},
	year={2009},
	publisher={SIAM}
}

@article{doerry2016bp,
	title={Basics of backprojection algorithm for processing synthetic aperture radar images},
	author={Doerry, Armin W and Bishop, Edward E and Miller, John A},
	journal={Sandia Report SAND2016-1682, Unlimited Release},
	pages={59},
	year={2016}
}

@article{feng2015synthetic,
	title={Synthetic aperture radar image despeckling via total generalised variation approach},
	author={Feng, Wensen and Lei, Hong and Qiao, Hong},
	journal={IET Image Processing},
	volume={9},
	number={3},
	pages={236--248},
	year={2015},
	publisher={Wiley Online Library}
}

@phdthesis{richardson2021multi,
	title={Multi-Bang Regularization and Applications},
	author={Richardson, Philip},
	year={2021},
	publisher={The University of Manchester (United Kingdom)}
}

@article{papafitsoros2015study,
	title={A study of the one dimensional total generalised variation regularisation problem},
	author={Papafitsoros, Konstantinos and Bredies, Kristian},
	journal={Inverse Problems and Imaging},
	volume={9},
	number={2},
	pages={511--550},
	year={2015},
	publisher={Inverse Problems and Imaging}
}

@article{jorgensen2021core1,
	title={Core Imaging Library-Part I: a versatile Python framework for tomographic imaging},
	author={J{\o}rgensen, Jakob S and Ametova, Evelina and Burca, Genoveva and Fardell, Gemma and Papoutsellis, Evangelos and Pasca, Edoardo and Thielemans, Kris and Turner, Martin and Warr, Ryan and Lionheart, William RB and others},
	journal={Philosophical Transactions of the Royal Society A},
	volume={379},
	number={2204},
	pages={20200192},
	year={2021},
	publisher={The Royal Society Publishing}
}

@article{papoutsellis2021core2,
	title={Core Imaging Library-Part II: multichannel reconstruction for dynamic and spectral tomography},
	author={Papoutsellis, Evangelos and Ametova, Evelina and Delplancke, Claire and Fardell, Gemma and J{\o}rgensen, Jakob S and Pasca, Edoardo and Turner, Martin and Warr, Ryan and Lionheart, William RB and Withers, Philip J},
	journal={Philosophical Transactions of the Royal Society A},
	volume={379},
	number={2204},
	pages={20200193},
	year={2021},
	publisher={The Royal Society Publishing}
}

@article{watson2022focusing,
	title={Focusing dynamic single-channel synthetic aperture radar video with optical flow-informed reconstruction},
	author={Watson, FM},
	journal={Electronics Letters},
	volume={58},
	number={25},
	pages={991--994},
	year={2022},
	publisher={Wiley Online Library}
}

@misc{watson2024resolving,
	title={Resolving Full-Wave Through-Wall Transmission Effects in Multi-Static Synthetic Aperture Radar}, 
	author={Francis Watson and Daniel Andre and William Robert Breckon Lionheart},
	year={2024},
	eprint={2403.10354},
	archivePrefix={arXiv},
	primaryClass={math.NA},
	note={to appear in IOP Inverse Problems}
}

@article{lin2001strongest,
	title={On strongest necessary and weakest sufficient conditions},
	author={Lin, Fangzhen},
	journal={Artificial Intelligence},
	volume={128},
	number={1-2},
	pages={143--159},
	year={2001},
	publisher={Elsevier}
}

@article{parikh2014proximal,
	title={Proximal algorithms},
	author={Parikh, Neal and Boyd, Stephen and others},
	journal={Foundations and trends{\textregistered} in Optimization},
	volume={1},
	number={3},
	pages={127--239},
	year={2014},
	publisher={Now Publishers, Inc.}
}

@misc{umbra,
	title={Umbra Synthetic Aperture Radar (SAR) Open Data was accessed on DATE from https://registry.opendata.aws/umbra-open-data, licensed under CC-BY-4.0}
}

@inproceedings{casteel2007challenge,
	title={A challenge problem for 2D/3D imaging of targets from a volumetric data set in an urban environment},
	author={Casteel Jr, Curtis H and Gorham, LeRoy A and Minardi, Michael J and Scarborough, Steven M and Naidu, Kiranmai D and Majumder, Uttam K},
	booktitle={Algorithms for Synthetic Aperture Radar Imagery XIV},
	volume={6568},
	pages={97--103},
	year={2007},
	organization={SPIE}
}

@article{beck2009siam,
	title={A fast iterative shrinkage-thresholding algorithm for linear inverse problems},
	author={Beck, Amir and Teboulle, Marc},
	journal={SIAM journal on imaging sciences},
	volume={2},
	number={1},
	pages={183--202},
	year={2009},
	publisher={SIAM}
}

@article{beck2009ieee,
	title={Fast gradient-based algorithms for constrained total variation image denoising and deblurring problems},
	author={Beck, Amir and Teboulle, Marc},
	journal={IEEE transactions on image processing},
	volume={18},
	number={11},
	pages={2419--2434},
	year={2009},
	publisher={IEEE}
}

@article{duran2016collaborative,
	title={Collaborative total variation: A general framework for vectorial TV models},
	author={Duran, Joan and Moeller, Michael and Sbert, Catalina and Cremers, Daniel},
	journal={SIAM Journal on Imaging Sciences},
	volume={9},
	number={1},
	pages={116--151},
	year={2016},
	publisher={SIAM}
}

@article{chambolle2011first,
	title={A first-order primal-dual algorithm for convex problems with applications to imaging},
	author={Chambolle, Antonin and Pock, Thomas},
	journal={Journal of mathematical imaging and vision},
	volume={40},
	pages={120--145},
	year={2011},
	publisher={Springer}
}

@inproceedings{stevens2017bright,
	title={Bright spark: Ka-band SAR technology demonstrator},
	author={Stevens, M and Jones, O and Moyse, P and Tu, S and Wilshire, A},
	booktitle={International Conference on Radar Systems (Radar 2017)},
	pages={1--4},
	year={2017},
	organization={IET}
}

@inproceedings{doody2017bright,
	author={Doody, S G and Hughes, N and Ramio-Tomas, L and Mak, E},
	booktitle={International Conference on Radar Systems (Radar 2017)}, 
	title={Bright sapphire - LF SAR imagery - A first look}, 
	year={2017},
	volume={},
	number={},
	pages={1-4},
	doi={10.1049/cp.2017.0461}
}

@article{gens1996review,
	title={Review Article SAR interferometry—issues, techniques, applications},
	author={Gens, Rudiger and Van Genderen, John L},
	journal={International journal of remote sensing},
	volume={17},
	number={10},
	pages={1803--1835},
	year={1996},
	publisher={Taylor \& Francis}
}

@article{rocca2000sar,
	title={SAR interferometry and its applications},
	author={Rocca, Fabio and Prati, Claudio and Monti Guarnieri, Andrea and Ferretti, Alessandro},
	journal={Surveys in Geophysics},
	volume={21},
	pages={159--176},
	year={2000},
	publisher={Springer}
}

@book{jakowatz2012spotlight,
	title={Spotlight-mode synthetic aperture radar: a signal processing approach: a signal processing approach},
	author={Jakowatz, Charles VJ and Wahl, Daniel E and Eichel, Paul H and Ghiglia, Dennis C and Thompson, Paul A},
	year={2012},
	publisher={Springer Science \& Business Media}
}

@inproceedings{andre2024moving,
	title={Moving Target Detection in Coherent Clutter with Reverse-Path Multistatic SAR},
	author={Andre, Daniel and Watson, Francis and Finnis, Mark},
	booktitle={EUSAR 2024; 15th European Conference on Synthetic Aperture Radar},
	pages={623--628},
	year={2024},
	organization={VDE}
}

\end{document}